\documentclass[12pt]{article}
\usepackage{exscale,relsize}
\usepackage[sort,nocompress]{cite}
\usepackage{amsmath}
\usepackage{amsfonts}
\usepackage[hidelinks]{hyperref}
\usepackage{amssymb}
\usepackage{calc}
\usepackage{theorem}
\usepackage{pifont}      
\usepackage{array}
\usepackage{color}

\usepackage{graphicx} 
\usepackage{float} 
\usepackage{subcaption} 

\newcommand{\N}{\mathbb{N}}   

\newcommand{\dist}{\ensuremath{\operatorname{dist}}} 

\newcommand{\norm}[1]{\left\lVert#1\right\rVert} 
\newcommand{\ip}[2]{\langle#1,#2\rangle} 

\newcommand{\gph}{\ensuremath{\operatorname{gph}}}

\newcommand{\R}{\ensuremath{\mathbb R}}
\newcommand{\Rn}{\ensuremath{\mathbb R^n}}

\newcommand{\dom}{\ensuremath{\operatorname{dom}}}

\newcommand{\inte}{\ensuremath{\operatorname{int}}}

\providecommand{\BB}[2]{\mathbb{B}(#1;#2)}

\newcommand{\prox}{\ensuremath{\operatorname{Prox}}}

\newtheorem{theorem}{Theorem}[section]
\newtheorem{lemma}[theorem]{Lemma}
\newtheorem{fact}[theorem]{Fact}
\newtheorem{corollary}[theorem]{Corollary}
\newtheorem{proposition}[theorem]{Proposition}
\newtheorem{defn}[theorem]{Definition}

\theoremstyle{plain}{\theorembodyfont{\rmfamily}
}
\theoremstyle{plain}{\theorembodyfont{\rmfamily}
}
\theoremstyle{plain}{\theorembodyfont{\rmfamily}
}
\theoremstyle{plain}{\theorembodyfont{\rmfamily}
\newtheorem{example}[theorem]{Example}}
\theoremstyle{plain}{\theorembodyfont{\rmfamily}
\newtheorem{remark}[theorem]{Remark}}
\theoremstyle{plain}{\theorembodyfont{\rmfamily}
}
\def\proof{\noindent{\it Proof}. \ignorespaces}
\def\endproof{\ensuremath{\quad \hfill \blacksquare}}

\newcommand{\pluss}{{\hskip1pt \raise1pt\vbox{\hrule width6pt \vskip1pt
\hrule width6pt}\kern-4pt{\lower1pt\hbox{\vrule height6pt \kern1pt\vrule
height6pt}}\hskip5pt}}

\newcommand{\argmin}{\mathop{\rm argmin}\limits}

\begin{document}
	\title{{\fontfamily{ptm}\selectfont The exact modulus of the generalized concave Kurdyka-\L ojasiewicz property}}

\author{
         Xianfu\ Wang\thanks{
                 Mathematics, University of British Columbia, Kelowna, B.C.\ V1V~1V7, Canada.
                 E-mail: \href{mailto:shawn.wang@ubc.ca}{\texttt{shawn.wang@ubc.ca}}.} and
         Ziyuan Wang\thanks{
                 Mathematics, University of British Columbia, Kelowna, B.C.\ V1V~1V7, Canada.
                 E-mail: \href{mailto:ziyuan.wang@alumni.ubc.ca}{\texttt{ziyuan.wang@alumni.ubc.ca}}.}
                 }


\maketitle

\begin{abstract} \noindent We introduce a generalized version of the concave Kurdyka-\L ojasiewicz (KL) property by employing nonsmooth desingularizing functions. We also present the exact modulus of the generalized concave KL property, which provides an answer to the open question regarding the optimal concave desingularizing function. The exact modulus is designed to be the smallest among all possible concave desingularizing functions. Examples are
given to illustrate this pleasant property. In turn, using the exact modulus
we provide the sharpest upper bound for the total length of iterates generated by
the celebrated Bolte-Sabach-Teboulle PALM algorithm.
\end{abstract}




\section{Introduction}\label{Intro}

The continuous optimization community has witnessed a surging interest of employing the concave KL property (see Definition~\ref{Def:KL property}) to solve problems from various applications, such as image processing~\cite{Ipiano2016,Banert2019}, compressed sensing~\cite{TKP2017extra,yu2020convergence,TKP2019}, machine learning~\cite{Lange2019} and many more. The aforementioned work, despite devoting to different proximal-type algorithms, share a common theme: Employing the concave KL property as a regularity condition to ensure the algorithm of interest has the finite length property; see, e.g.,~\cite[Theorem 4.9]{Ipiano2016} and \cite[Theorem 3.1]{Lange2019}. This pleasant convergence methodology can be traced back to the fundamental work of Bolte et~al.~\cite{bolte2010survey,bolte2007lojasiewicz} and Attouch et~al.~\cite{Attouch2010}.

In the concave KL property, the
 concave desingularizing function plays a central role in
 estimating both the convergence rate and total length of iterates
 generated by the algorithm of interest; see, e.g.,~\cite[Theorem 1, Lemma 4]{Banert2019}. However, the concave desingularizing functions
  are not necessarily unique. It is natural to ask what the optimal~(minimal) one is. This question remains open in the current literature. Classic definition of the concave KL property requires continuous differentiability
 of desingularizing functions, precluding the infimum of all concave desingularizing functions from staying within the same class.
This paper is devoted to answering the open question:
\begin{equation}\label{open question}
	\emph{What is the optimal concave desingularizing function?}	
\end{equation}
To this end, we introduce an extension of the concave KL property and its associated exact modulus by allowing nonsmooth desingularizing functions. This extended framework allows us to capture the optimal concave desingularizing function
through the exact modulus, yet still compatible with the usual concave KL convergence technique employed by a vast amount of literature. Our work opens the door to improve
 convergence results of a broad range of algorithms that adopt the concave KL assumption.

Throughout this paper, $\Rn$ is the standard Euclidean space
with inner product $\ip{x}{y}=x^Ty$ and Euclidean norm $\norm{x}=\sqrt{\ip{x}{x}}$ for $x, y\in\Rn$.
The open ball centered at $\bar{x}\in \Rn$ with radius $r>0$ is denoted by~$\BB{\bar{x}}{r}$.
We let $\overline{\R}=(-\infty,\infty]$, $\R_+=[0,\infty)$, and $\N=\{1,2,3,\ldots\}$.
The distance function of a subset $K\subseteq \Rn$ is $\dist(\cdot,K):\Rn\rightarrow [0,\infty]:
x\mapsto\dist(x,K)=\inf\{\norm{x-y}:y\in K\}.$
For $f:\Rn\to\overline{\R}$ and $r_1,r_2\in[-\infty,\infty]$, we set $[r_1<f<r_2]=\{x\in\Rn:r_1<f(x)<r_2\}$. For $\eta\in(0,\infty]$, denote by $\mathcal{K}_\eta$ the class of functions $\varphi:[0,\eta)\rightarrow\R_+$ that satisfy the following three conditions: (i) $\varphi:[0,\eta)\rightarrow\R_+$ is continuous with $\varphi(0)=0$; (ii) $\varphi$ is $C^1$ on $(0,\eta)$; (iii) $\varphi^\prime(t)>0$ for all $t\in(0,\eta)$. The pointwise version\footnote{In the remainder of this paper, we shall simply refer this pointwise definition as ``the concave KL property" for the sake of simplicity. However we would like to remind readers that the concave KL property is originally introduced as a property about function values instead of points; see~\cite[Theorem 14]{bolte2007}.} of the concave KL property is defined as follows.
\begin{defn}\label{Def:KL property}
	Let $f:\mathbb{R}^n\rightarrow\overline{\mathbb{R}}$ be proper and lower semicontinuous (lsc).
	
	(i) We say $f$ has the KL property at $\bar{x}\in\dom\partial f$ if there exist a
 neighborhood $U\ni\bar{x}$, $\eta\in(0,\infty]$ and a function $\varphi\in\mathcal{K}_\eta$ such that for all $x\in U\cap[0<f-f(\bar{x})<\eta]$,
	\begin{equation}\label{Inequality: KL inequality}
	\varphi^\prime\big(f(x)-f(\bar{x})\big)\cdot\dist\big(0,\partial f(x)\big)\geq1,
	\end{equation}		
where $\partial f(x)$ denotes the limiting subdifferential of $f$ at $x$ (see Definition~\ref{Defn:limiting subdifferential}). The function $\varphi$ is called a \textit{desingularizing function} of $f$ at $\bar{x}$ with respect to $U$ and $\eta$. We say $f$ is a KL function if it has the KL property at every $\bar{x}\in\dom\partial f$.

(ii) We say that $f$ has the concave KL property at $\bar{x}\in\dom\partial f$ if it has the KL property at $\bar{x}$ with desingularizing function $\varphi\in\mathcal{K}_\eta$ being concave. Moreover, we say $f$ is a concave KL function if it has the concave KL property at every $\bar{x}\in\dom\partial f$.
\end{defn}

The pioneering work of \L ojasiewicz~\cite{Lojas1963} and Kurdyka~\cite{Kur98} on
 differentiable functions laid the foundation of the KL property, which was extended to nonsmooth functions by Bolte et al.\ in \cite{bolte2007lojasiewicz,bolte2007}.
 In the seminal work~\cite{bolte2010survey}, Bolte et al.\ coined the term ``KL property", gave characterizations and proposed the BDLM desingularizing function, which is the optimal desingularizing function
 under certain continuity and locally integrability conditions; see Fact~\ref{fact: integral condition} for details and \cite[Theorem 1]{Kur98} for a similar result in a different setting. However, the optimal concave desingularizing function associated with the concave KL property
 may not be captured by the BDLM desingularizing function when the continuity and integrability assumptions fail; see Section~\ref{sec: Comparison to integrability condition }.

The main contributions of this paper are listed below:
\begin{itemize}
\item Definition~\ref{Def: g-KL} generalizes the concave KL property. The main difference is that we allow the desingularizing function to be {non-differentiable}.
\item Proposition~\ref{Prop: g-KL modulus is optimal} shows that the exact modulus of the generalized concave KL property, given in Definition~\ref{Defn: g-KL modulus}, is the optimal concave desingularizing function, provided the existence of concave desingularizing functions. This result answers the open question~(\ref{open question}). 
\item Theorem~\ref{Theorem: finite length property} provides the sharpest upper bound on $\sum_{k=1}^\infty\norm{z_{k+1}-z_k}$, where $(z_k)_{k\in\N}$ is a sequence generated by the PALM algorithm. This result improves \cite[Theorem 1]{Bolte2014}.
\end{itemize}
\noindent Although most published articles emphasize desingularizing functions of the form $\varphi(t)=c\cdot t^{1-\theta}$ for $c>0$ and $\theta\in[0,1)$, the exact modulus has various forms.
Proposition~\ref{Prop: modulus of C1 convex} gives an explicit formula for the
optimal concave desingularizing function of locally convex and $C^1$ functions on the real line, in which case the exact modulus coincides with the desingularizing function obtained from the BDLM integrability condition. However, examples are given to show that the exact modulus is indeed the smaller one, even for nondifferentiable convex functions on the real line; see Examples~\ref{ex: campare to integral condition with convexity} and \ref{ex:piewiselinear example}. More examples comparing these two objects are provided in Section~\ref{sec: Comparison to integrability condition }. 
As a by-product concerning intersections of convex functions, we show in Example~\ref{Prop: inf of intersection} that there exist distinct strictly increasing convex $C^2$ functions $f, g:\R_+\rightarrow\R_+$
with $f(0)=g(0)=0$ such that $\inf\{x>0:f(x)=g(x)\}=0$.
Using our technique in Theorem~\ref{Theorem: finite length property}, one may improve other algorithms that adopt the concave KL property assumption.

The structure of this paper is as the following: Elements in variational analysis, classical analysis and facts of the classical KL property are collected in Section~\ref{Section: Auxiliary results}. The generalized concave KL property, the exact modulus and their properties are studied in Section~\ref{Section: g-KL and exact modulus}. Various examples and comparisons to the BDLM desingularizing functions are also given in that section. We revisit the celebrated PALM algorithm in Section~\ref{Section: PALM}. Concluding remarks and directions for future work are presented in Section~\ref{Sec:conclusion}.

\section{Preliminaries}\label{Section: Auxiliary results}
\subsection{Elements of variational and classical analysis}

We will use frequently the following subgradients in the nonconvex setting; see, e.g., \cite{rockwets, mor2005variational}.

\begin{defn}\label{Defn:limiting subdifferential}
	Let $f:\Rn\rightarrow\overline{\R}$ be a proper function. We say that
	\begin{itemize}
		\item[(i)] $v\in\Rn$ is a \textit{Fr\'echet subgradient} of $f$ at $\bar{x}\in\dom f$, denoted by $v\in\hat{\partial}f(\bar{x})$, if for every $x\in\dom f$,
		\begin{equation}\label{Formula:frechet subgradient inequality}
		f(x)\geq f(\bar{x})+\ip{v}{x-\bar{x}}+o(\norm{x-\bar{x}}).
		\end{equation}
		\item[(ii)] $v\in\Rn$ is a \textit{limiting subgradient} of $f$ at $\bar{x}\in\dom f$, denoted by $v\in\partial f(\bar{x})$, if
		\begin{equation}\label{Formula:limiting subgraident definition}
		v\in\{v\in\Rn:\exists x_k\xrightarrow[]{f}\bar{x},\exists v_k\in\hat{\partial}f(x_k),v_k\rightarrow v\},
		\end{equation}
where $x_k\xrightarrow[]{f}\bar{x}\Leftrightarrow x_k\rightarrow\bar{x}\text{ and }f(x_k)\rightarrow f(\bar{x})$. Moreover, we set $\dom\partial f=\{x\in\Rn:\partial f(x)\neq\emptyset\}$. We say that $\bar{x}\in\dom\partial f$ is a stationary point if $0\in\partial f(\bar{x})$.
	\end{itemize}
\end{defn}

Below is the definition of proximal mapping; see, e.g.,~\cite{Bolte2014} and~\cite{rockwets}.
\begin{defn}\label{Def:Moreau envelope and prox mapping} Let $f:\Rn\rightarrow\overline{\R}$ be proper and lsc and let $\lambda$ be a positive real number. The proximal mapping is defined by
\begin{align*}
(\forall x\in\Rn)\ \prox_\lambda^f(x)&=\argmin_{y\in\Rn}\Big\{f(y)+\frac{\lambda}{2}\norm{x-y}^2\Big\}.
\end{align*}
\end{defn}
The following fact follows from \cite[Theorem 1.25]{rockwets}.
\begin{fact}\label{fact:prox is well defined}  Let $f:\Rn\rightarrow\overline{\R}$ be proper and lsc with $\inf_{\Rn}f>-\infty$. Then for $\lambda\in(0,\infty)$, $\prox_\lambda^f(x)$ is nonempty for every $x\in\Rn$. Moreover, for every $v\in\Rn$ and $x\in\Rn$ we have
\[\prox_\lambda^f\left(x-\frac{1}{\lambda}v\right)=\argmin_y\Big\{\ip{y-x}{v}+\frac{\lambda}{2}\norm{x-y}^2+f(y)\Big\}.\]
\end{fact}

Some well-known properties of convex functions on the real line are given in the following fact.
\begin{fact}\label{Lemma: properties of convex functions on the line}\emph{(\cite[Section 24]{Rockafellar},~\cite[Chapter 17]{BC})} Let $I\subseteq\R$ be an open interval and let $\varphi:I\rightarrow\R$ be convex. Then

(i)  The side derivatives $\varphi_-^\prime(t)$ and $\varphi_+^\prime(t)$ are finite at every $t\in I$. Moreover, $\varphi_-^\prime(t)$ and $\varphi_+^\prime(t)$ are increasing.

(ii) $\varphi$ is differentiable except at countably many points of $I$, and $\varphi(s)-\varphi(t)=
\int_{t}^{s}\varphi'_{-}(x)dx=\int_{t}^{s}\varphi'_{+}(x)dx$
for all $s,t\in I$.

(iii) Let $t\in I$. Then for every $s\in I$, $\varphi(s)-\varphi(t)\geq \varphi_-^\prime(t)\cdot(s-t)$.
\end{fact}

The following result concerns the absolute continuity of integrals.
\begin{fact}\label{Lemma: Lebesgue integral vanishes}\emph{\cite[Theorem 6.79]{Stromberg}} Let $f\in L^1$. Then for each $\varepsilon>0$, there exists $\delta>0$ such that \[\int_{E}|f(x)|ds<\varepsilon,\]
	whenever $m(E)<\delta$, where $E$ is a Lebesgue-measurable set and $m(E)$ denotes its Lebesgue measure.
\end{fact}

\subsection{The Kurdyka-\L ojasiewicz property and known desingularizing functions}
In this section we collect several facts about the KL property and desingularizing functions. We begin with a result asserting that the KL property at non-stationary points is automatic; see, e.g.,~\cite[Remark 3.2(b)]{Attouch2010} and also~\cite[Lemma 2.1]{li2018calculus} for a detailed proof.

\begin{fact}\label{Fact: KL at non-stationary point} Let $f:\mathbb{R}^n\rightarrow\overline{\R}$ be proper and lsc. Let $\bar{x}\in\dom\partial f$ be a non-stationary point. Then there exist $c>0$ and $\theta\in[0,1)$ such that $f$ has the KL property at $\bar{x}$ with respect to $U=\BB{\bar{x}}{\varepsilon}$, $\eta=\varepsilon$ and $\varphi(t)=c\cdot t^{1-\theta}$.
\end{fact}

We now recall some desingularizing functions described by Bolte et~al.~\cite{bolte2010survey}, which we will later compare with our main results. Recall that a proper and lsc function $f:\Rn\to\overline{\R}$ is semiconvex if there exists $\alpha>0$ such that $f+\frac{\alpha}{2}\norm{\cdot}^2$ is convex. 

\begin{fact}\emph{\cite[Lemma 45, Theorem 18]{bolte2010survey}}\label{fact: integral condition} Let $f:\Rn\to\overline{\R}$ be lsc and semiconvex. Let $\bar{x}\in[f=0]$ and assume that there exist $\bar{r},\bar{\varepsilon}>0$ such that
	\begin{align}\label{formula:noncriticality assumption}
		x\in\BB{\bar{x}}{\bar{\varepsilon}}\cap[0<f\leq\bar{r}]\Rightarrow	0\notin\partial f(x).
	\end{align}
Suppose there exist $r_0\in(0,\bar{r})$ and $\varepsilon\in(0,\bar{\varepsilon})$ such that the function
\begin{equation}\label{formula: integral condition}
	u(r)=\frac{1}{\displaystyle\inf_{x\in\BB{\bar{x}}{\varepsilon}\cap[f=r]}\dist(0,\partial f(x))},~r\in(0,r_0]
\end{equation}
is finite-valued and belongs to $L^1(0,r_0)$. Then the following statements hold:

(i) There exists a continuous majorant $\bar{u}:(0,r_0]\to(0,\infty)$ such that $\bar{u}\in L^1(0,r_0)$ and $\bar{u}(r)\geq u(r)$ for all $r\in(0,r_0]$.
	
(ii)  Define for $t\in(0,\bar{r})$ \[\varphi(t)=\int_0^t\bar{u}(s)ds.\]
Then $\varphi\in\mathcal{K}_{r_0}$. For every $x\in\BB{\bar{x}}{\varepsilon}\cap[0<f\leq r_0]$, one has
\[\varphi^\prime\left(f(x)\right)\dist(0,\partial f(x))\geq1.\]

	
\end{fact}
\begin{remark}\label{rem: limitation of integral condition} Fact~\ref{fact: integral condition} is extracted from the implication $(v)\Rightarrow(i)$ in the proof of~\cite[Theorem 18]{bolte2010survey}, where the above desingularizing function $\varphi(t)$ was not stated explicitly in their theorem statement. Since results in this paper are on $\Rn$, we restrict Fact~\ref{fact: integral condition} to $\Rn$, in which case Assumption (24) of \cite[Theorem 18]{bolte2010survey} becomes superfluous.
\end{remark}

Next we collect facts that ensure existence of concave desingularizing functions, which set the stage for our main results. With convexity, the following fact asserts that the desingularizing function given by Fact~\ref{fact: integral condition} can be taken to be concave with an enlarged domain $[0,\infty)$.

\begin{fact}\emph{\cite[Lemma 45, Theorem 29]{bolte2010survey}}\label{fact:integral condition with convexity} Let $f:\Rn\rightarrow\overline{\R}$ be a proper, lsc and convex function with $\inf f=0$. Suppose that there exist $r_0>0$ and $\varphi\in\mathcal{K}_{r_0}$ such that for all $x\in[0<f\leq r_0]$,
\[\varphi^\prime(f(x))\dist(0,\partial f(x))\geq1.\]
Then the following statements hold:

(i) Define for $r\in(0,\infty)$ the function \[u(r)=\frac{1}{\displaystyle\inf_{x\in[f=r]}\dist(0,\partial f(x))}.\]
Then $u$ is finite-valued, decreasing and $u\in L^1(0,r_0)$. Moreover, there exists a decreasing continuous function $\tilde{u}\in L^1(0,r_0)$ such that $\tilde{u}\geq u$.

(ii) Pick $\bar{r}\in(0,r_0)$ and define for $r\in(0,\infty)$
\begin{align*}
 \varphi(r)=
\begin{cases}
\int_0^r\tilde{u}(s)ds,&\text{if }r\leq \bar{r};\\
\int_0^{\bar{r}}\tilde{u}(s)ds+\tilde{u}(\bar{r})(r-\bar{r}),& \text{otherwise}.
\end{cases}
\end{align*}
Then $\varphi\in\mathcal{K}_\infty$ is concave and for every $x\notin[f=0]$,
\[ \varphi^\prime(f(x))\dist(0,\partial f(x))\geq1.\]
\end{fact}

Another celebrated result states that semialgebraic functions have the concave KL property.
\begin{defn} (i) A set $E\subseteq\mathbb{R}^n$ is called \textit{semialgebraic} if there exist finitely many polynomials $g_{ij}, h_{ij}:\mathbb{R}^n\rightarrow\mathbb{R}$ such that	
\begin{equation*}		E=\bigcup_{j=1}^p\bigcap_{i=1}^q\{x\in\mathbb{			R}^n:g_{ij}(x)=0\text{ and }h_{ij}(x)<0\}.	
\end{equation*}	

(ii) A function $f:\mathbb{R}^n\rightarrow\overline{\R}$ is called \textit{semialgebraic} if its graph	
\begin{equation*}		\gph f=\{(x,y)\in\mathbb{R}^{n+1}:f(x)=y\}	\end{equation*}
is semialgebraic.
\end{defn}

\begin{fact}\emph{\cite[Corollary 16]{bolte2007}}\label{Fact:semi-algebraic functions are KL} Let $f:\mathbb{R}^n\rightarrow\overline{\R}$ be a proper and lsc function and let $\bar{x}\in\dom\partial f$. If $f$ is semialgebraic, then it has the concave KL property at $\bar{x}$ with $\varphi(t)=c\cdot t^{1-\theta}$ for some $c>0$ and $\theta\in(0,1)$.
\end{fact}
\begin{remark} (i) Many useful functions in optimization are semialgebraic; see,~e.g.,~\cite{Attouch2010,Bolte2014} and the references therein. Functions definable in o-minimal structure, which include semialgebraic functions, also satisfy the concave KL property; see
~\cite{Attouch2010,bolte2007}.

(ii) Although it is well-known that real-polynomials are semialgebraic and thus have the KL property, only until very recently
\cite[Corollary 9]{bolte2016error} did Bolte et al.\  provide an explicit formula for desingularizing functions of convex piecewise polynomials.
\end{remark}

One can also determine the desingularizing function of the KL property for convex functions through the following growth condition.
\begin{fact}\emph{\cite[Theorem 30]{bolte2010survey}}\label{Fact:Growth condition gives KL} Let $f:\Rn\to\overline{\R}$ be a proper lsc convex function with $f(0)=\min f$. Let $S\subseteq\Rn$. Assume that there exists a function $m:\R_+\to\R_+$ that is continuous, strictly increasing, $m(0)=0$, $f\geq m\left(\dist(\cdot,\argmin f) \right)$ on $S\cap\dom f$ and  \[\exists\rho>0, ~\int_0^\rho\frac{m^{-1}(s)}{s}ds<\infty,\]
where $m^{-1}$ denotes the inverse function of $m$. Then for all $x\in S\backslash\argmin f$,
\[\varphi^\prime(f(x))\dist\left(0,\partial f(x)\right)\geq1,\]
where for $t\in(0,\rho)$,
$$\varphi(t)=\int_0^t\frac{m^{-1}(s)}{s}ds.$$
\end{fact}

\begin{remark} Assuming there exists a concave desingularizing function, Facts~\ref{fact: integral condition},~\ref{fact:integral condition with convexity} and~\ref{Fact:Growth condition gives KL} may fail to capture the optimal one, even for convex functions on the real line; see Section~\ref{sec: Comparison to integrability condition }.
\end{remark}

\section{The generalized concave KL property and its exact modulus}\label{Section: g-KL and exact modulus}
In this section, we provide an answer to the open question~(\ref{open question}). The generalized concave Kurdyka-\L ojasiewicz property and its exact modulus are introduced to provide the answer. Given existence of concave desingularizing functions, we shall see that the exact modulus is indeed optimal.
\subsection{The generalized concave KL property}
For $\eta\in(0,\infty]$, denote by $\Phi_\eta$ the class of functions $\varphi:[0,\eta)\rightarrow\R_+$ satisfying the following conditions: (i) $\varphi(t)$ is right-continuous at $t=0$ with $\varphi(0)=0$; (ii) $\varphi$ is strictly increasing on $[0,\eta)$. Recall that the left derivative of $\varphi:[0,\infty)\rightarrow\R$ at $t\in(0,\infty)$ is defined by
\[\varphi_-^\prime(t)=\lim_{s\rightarrow t^-}\frac{\varphi(s)-\varphi(t)}{s-t}.\]

Some useful properties of concave $\varphi\in\Phi_\eta$ are collected below.

\begin{lemma}\label{left derivative inequality} For $\eta\in(0,\infty]$ and concave $\varphi\in\Phi_\eta$, the following assertions hold:
	\begin{itemize}
		\item[(i)] Let $t>0$. Then $\varphi(t)=\lim_{u\rightarrow0^+}\int_u^t\varphi_-^\prime(s)ds=\int_0^t\varphi_-^\prime(s)ds$.
		\item [(ii)] The function $t\mapsto\varphi_-^\prime(t)$ is decreasing and $\varphi_-^\prime(t)>0$ for $t\in(0,\eta)$. 
		\item[(iii)] For $0\leq s<t<\eta$, $\varphi^\prime_-(t)\leq\frac{\varphi(t)-\varphi(s)}{t-s}$.
	\end{itemize}
\end{lemma}
\proof (i) Invoking Fact~\ref{Lemma: properties of convex functions on the line}(ii) yields
\[\varphi(t)=\lim_{u\rightarrow0^+}\big(\varphi(t)-\varphi(u)\big)=\lim_{u\rightarrow0^+}\int_{u}^t\varphi_-^\prime(s)ds<\infty,\]
where the first equality holds because $\varphi$ is right-continuous at $0$ with $\varphi(0)=0$. Let $(u_n)_{n\in\N}$ be a decreasing sequence with $u_1<t$ such that $u_n\rightarrow0^+$ as $n\rightarrow\infty$. For each $n$, define $h_n:(0,t]\rightarrow\R_+$ by $h_n(s)=\varphi_-^\prime(s)$ if $s\in(u_n,t]$ and $h_n(s)=0$ otherwise. Then the sequence $(h_n)_{n\in\N}$ satisfies: (a) $h_n\leq h_{n+1}$ for every $n\in\N$; (b) $h_n(s)\rightarrow\varphi_-^\prime(s)$ pointwise on $(0,t)$; (c) The integral $\int_0^th_n(s)ds=\int_{u_n}^t\varphi_-^\prime(s)ds=\varphi(t)-\varphi(u_n)\leq \varphi(t)-\varphi(0)<\infty$ for every $n\in\N$. Hence the monotone convergence theorem implies that
\[\lim_{u\rightarrow0^+}\int_{u}^t\varphi_-^\prime(s)ds=\lim_{n\rightarrow\infty}\int_{u_n}^t\varphi_-^\prime(s)ds=\lim_{n\rightarrow\infty}\int_0^th_n(s)ds=\int_0^t\varphi_-^\prime(s)ds.\]

(ii) According to Fact~\ref{Lemma: properties of convex functions on the line}(i), the function $t\mapsto\varphi_-^\prime(t)$ is decreasing. Suppose that $\varphi_-^\prime(t_0)=0$ for some $t_0\in(0,\eta)$. Then by the monotonicity of $\varphi_-^\prime$ and (i), we would have $\varphi(t)-\varphi(t_0)=\int_{t_0}^t\varphi_-^\prime(s)ds\leq(t-t_0)\varphi_-^\prime(t_0)=0$ for $t>t_0$, which contradicts to the assumption that $\varphi$ is strictly increasing.

(iii) For $0<s<t<\eta$, applying Fact~\ref{Lemma: properties of convex functions on the line}(iii) to the convex function $-\varphi$ yields that $-\varphi(s)+\varphi(t)\geq-\varphi_-^\prime(t)(s-t)\Leftrightarrow\varphi_-^\prime(t)\leq\big(\varphi(t)-\varphi(s)\big)/(t-s)$. The desired inequality then follows from the right-continuity of $\varphi$ at $0$.\endproof

Now we introduce the pointwise generalized concave KL property and its setwise variant.

\begin{defn}\label{Def: g-KL} Let $f:\Rn\rightarrow\overline{\mathbb{R}}$ be proper and lsc.  Let $\bar{x}\in\dom\partial f$ and $\mu\in\R$, and let $V\subseteq\dom\partial f$ be a nonempty subset.

(i) We say that $f$ has the pointwise \textit{generalized concave KL} at $\bar{x}\in\dom\partial f$ if there exist a neighborhood $U\ni\bar{x}$, $\eta\in(0,\infty]$ and concave $\varphi\in\Phi_\eta$, such that for all $x\in U\cap[0<f-f(\bar{x})<\eta]$,
	\begin{equation}\label{g-KL inequality}
	\varphi^\prime_-\big(f(x)-f(\bar{x})\big)\cdot\dist\big(0,\partial f(x)\big)\geq1.
	\end{equation}
	
(ii) Suppose that $f(x)=\mu$ on $V$. We say $f$ has the setwise generalized concave KL property\footnote{For simplicity, we shall omit adjectives ``pointwise" and ``setwise" whenever there is no ambiguity. We also remind readers that one can define similarly generalized concave KL property around function values, which, however, will not be treated in the paper.} on $V$ if there exist $U\supset V$, $\eta\in(0,\infty]$ and concave $\varphi\in\Phi_\eta$ such that for every $x\in U\cap[0<f-\mu<\eta]$,\begin{equation}\label{uniform g-KL inequality}	\varphi^\prime_-\big(f(x)-\mu\big)\cdot\dist\big(0,\partial f(x)\big)\geq1.\end{equation}
\end{defn}
\begin{remark} (i) Evidently the generalized concave KL property on a set $V$ reduces to the generalized concave KL property at $\bar{x}$ if $V=\{\bar{x}\}$. This setwise definition will be useful in Section~\ref{Section: PALM}.
	
(ii) Clearly, the concave KL property (see Definition~\ref{Def:KL property}) implies the generalized concave KL property. However, the generalized notion allows desingularizing functions to be non-differentiable by using the left derivative, which is well-defined thanks to Fact~\ref{Lemma: properties of convex functions on the line}.
\end{remark}

In the rest of this subsection, we work towards generalizing a result by Bolte et al.~\cite[Lemma 6]{Bolte2014}, whose proof we will follow. For nonempty subset $ \Omega\subseteq\Rn$ and $\varepsilon\in(0,\infty]$, define $ \Omega_\varepsilon=\{x\in\Rn:\dist(x, \Omega)<\varepsilon \}.$
Let us recall the Lebesgue number lemma~\cite[Theorem 55]{Pugh}.
\begin{lemma}\label{Lemma: Lebesgue number lemma}
	Let $ \Omega\subseteq\Rn$ be a nonempty compact subset. Suppose that $\{U_i\}_{i=1}^p$ is a finite open cover of $ \Omega$.
Then there exists $\varepsilon>0$, which is called the Lebesgue number of $\Omega$, such that \[ \Omega\subseteq  \Omega_\varepsilon\subseteq\bigcup_{i=1}^pU_i.\]
\end{lemma}

Proposition~\ref{Prop: Uniformize the g-KL} below connects the pointwise generalized concave KL property to its setwise counterpart, generalizes~\cite[Lemma 6]{Bolte2014}, and will play a key role in Section~\ref{Section: PALM}.

\begin{proposition}\label{Prop: Uniformize the g-KL} Let $f:\Rn\rightarrow\overline{\R}$ be proper lsc and let $\mu\in\R$. Let $\Omega\subseteq\dom\partial f$ be a nonempty compact set on which $f(x)=\mu$ for all $x\in\Omega$. Suppose that $f$ satisfies the pointwise generalized concave KL property at each $x\in\Omega$. Then there exist $\varepsilon>0,\eta\in(0,\infty]$ and concave $\varphi(t)\in\Phi_\eta$ such that $f$ has the setwise generalized concave KL property on $\Omega$ with respect to $U=\Omega_\varepsilon$, $\eta$ and $\varphi$.
\end{proposition}
\proof For each $x\in\Omega$, there exist $\varepsilon=\varepsilon(x)>0$, $\eta=\eta(x)\in(0,\infty]$ and concave $\varphi(t)=\varphi_x(t)\in\Phi_\eta$ such that for $y\in\BB{x}{\varepsilon}\cap[0<f-f(x)<\eta]$, \[\varphi_-^\prime\big(f(y)-f(x)\big)\cdot\dist\big(0,\partial f(y)\big)\geq1.\]
Note that $\Omega\subseteq\bigcup_{x\in\Omega}\BB{x}{\varepsilon}$. Because $\Omega$ is compact, there exist elements $x_1,\ldots,x_p\in\Omega$ such that $\Omega\subseteq\bigcup_{i=1}^p\BB{x_i}{\varepsilon_i}$. Moreover, for each $i$ and  $x\in\BB{x_i}{\varepsilon_i}\cap[0<f-f(x_i)<\eta_i]=\BB{x_i}{\varepsilon_i}\cap[0<f-\mu<\eta_i]$, one has
\begin{equation}\label{wwww}
(\varphi_i)_-^\prime\big(f(x)-\mu\big)\cdot\dist\big(0,\partial f(x)\big)\geq1.
\end{equation}

Define $\varphi(t)=\sum_{i=1}^p\varphi_i(t)$ and $\eta=\min_{1\leq i\leq p}\eta_i$. Evidently, $\varphi$ is a concave function belonging to $\Phi_\eta$. By
Lemma~\ref{Lemma: Lebesgue number lemma}, there exists $\varepsilon>0$ such that $\Omega\subseteq\Omega_\varepsilon\subseteq\bigcup_{i=1}^p\BB{x_i}{\varepsilon_i},$
which by the fact that $\eta\leq\eta_i$ for every $i$ further implies that
\[x\in\Omega_\varepsilon\cap[0<f-\mu<\eta]\Rightarrow\exists i_0,~s.t.,~x\in\BB{x_{i_0}}{\varepsilon_{i_0}}\cap[0<f-\mu<\eta_{i_0}].\]
Hence for every $x\in\Omega_\varepsilon\cap[0<f-\mu<\eta]$, one has
\begin{align*}
\varphi_-^\prime\big(f(x)-\mu\big)\cdot\dist\big(0,\partial f(x)\big)\geq(\varphi_{i_0})_-^\prime\big(f(x)-\mu\big)\cdot\dist\big(0,\partial f(x)\big)\geq1,
\end{align*}
where the first inequality holds because $(\varphi_i)_-^\prime(t)>0$ by Lemma~\ref{left derivative inequality}, and the second inequality is implied by (\ref{wwww}).\endproof

\subsection{The exact modulus of the generalized concave KL property}

Following the definition of the generalized concave KL property, we introduce its associated exact modulus.

\begin{defn}\label{Defn: g-KL modulus} Let $f:\Rn\rightarrow\overline{\R}$ be proper and lsc. Let $\bar{x}\in\dom\partial f$ and let $U\subseteq\dom\partial f$ be a neighborhood of $\bar{x}$. Let $\eta\in(0,\infty]$. Furthermore, define $h:(0,\eta)\rightarrow\R$ by
	\[h(s)=\sup\big\{\dist^{-1}\big(0,\partial f(x)\big): x\in U\cap[0<f-f(\bar{x})<\eta],s\leq f(x)-f(\bar{x})\big\}.\]
	Suppose that $h(s)<\infty$ for $s\in(0,\eta)$. The \emph{exact modulus of the generalized concave KL property} of $f$ at $\bar{x}$ with respect to $U$ and $\eta$ is the function
	$\tilde{\varphi}:[0,\eta)\rightarrow\R_+: t\mapsto\int_0^th(s)ds,~\forall t\in(0,\eta),$
 and $\tilde{\varphi}(0)=0$. If $U\cap[0<f-f(\bar{x})<\eta]=\emptyset$ for given $U\ni\bar{x}$ and $\eta>0$,
 then we set the exact modulus with respect to $U$ and $\eta$ to be $\tilde{\varphi}(t)\equiv0$.
\end{defn}
\begin{remark}\label{rem:differences compared to integral condition}(i) The essential difference between the exact modulus and the BDLM desingularizing function in Fact~\ref{fact: integral condition} is that the exact modulus utilizes the set $U\cap[s\leq f(x)-f(\bar{x})]$ instead of $U\cap[s=f(x)-f(\bar{x})]$. In addition, the exact modulus $\tilde{\varphi}$ is not necessarily differentiable while the BDLM desingularizing requires differentiability. In order for the exact modulus to be well-defined, however, it requires the existence of concave desingularizing functions,
which is a strong assumption. Examples of such functions include convex functions satisfying the KL property and semialgebraic functions, cf.\ Facts~\ref{fact:integral condition with convexity} and~\ref{Fact:semi-algebraic functions are KL}, which are frequently treated in papers devoted to algorithmic applications of the concave KL property~\cite{Attouch2010,Bolte2014,Lange2019,yu2020convergence,Ipiano2016,TKP2019,TKP2017extra,Banert2019}.
		
(ii) Note that $\lim_{s\rightarrow0^+}h(s)$ should be infinity if $\bar{x}$ is a stationary point, in which case the function $\tilde{\varphi}(t)$ represents a limit of Riemann or Lebesgue integrals.
	
(iii) The assumption that $h(s)<\infty$ for $s\in(0,\eta)$ is necessary. For example, consider the exact modulus of the generalized concave KL property of the function $f(x)=1-e^{-|x|}$ at $0$. Then one has for $x\neq0$, $\dist^{-1}\big(0,\partial f(x)\big)=e^{|x|}$. Let $U=\R$ and $\eta_1=1$. Then
	\begin{align*}
	h_1(s)&=h_{U,\eta_1}(s)=\sup\big\{\dist^{-1}\big(0,\partial f(x)\big):\R\cap[0<f<1],s\leq f(x)\big\}=\infty.
	\end{align*}
	This can be avoided by shrinking the set $U\cap[0<f-f(\bar{x})<\eta]$. Let $\eta_2\in(0,1)$. Then
	\begin{align*}
	h_2(s)&=h_{U,\eta_2}(s)=\sup\big\{\text{dist}^{-1}(0,\partial f(x)):x\in\mathbb{R}\cap[0<f<\eta_2],s\leq f(x)\big\}=\frac{1}{1-\eta_2}.
	\end{align*}
\end{remark}

The exact modulus $\tilde{\varphi}$ is designed to be the optimal concave desingularizing function. The following lemma is needed to prove this property.

\begin{lemma}\label{Lemma:integration of increasing function is convex}
Let $\eta\in(0,\infty]$ and let $h:(0,\eta)\rightarrow\R_+$ be a positive-valued decreasing function. Define $\varphi(t)=\int_0^th(s)ds$ for $t\in(0,\eta)$ and set $\varphi(0)=0$. Suppose that $\varphi(t)<\infty$ for $t\in(0,\eta)$. Then $\varphi$ is a strictly increasing concave function on $[0,\eta)$ with $$\varphi_-^\prime(t)\geq h(t)$$ for $t\in(0,\eta)$, and right-continuous at $0$. If in addition $h$ is a continuous function, then $\varphi$ is $C^1$ on $(0,\eta)$.
\end{lemma}
\proof Let $0<t_0<t_1<\eta$. Then $\varphi(t_1)-\varphi(t_0)=\int_{t_0}^{t_1}h(s)ds\geq(t_1-t_0)\cdot h(t_1)>0$, which means $\varphi$ is strictly increasing. Applying Fact~\ref{Lemma: Lebesgue integral vanishes}, one concludes that $\varphi(t)\rightarrow\varphi(0)=0$ as $t\rightarrow0^+$. The concavity of $\varphi$ and the inequality $\varphi_-^\prime(t)\geq h(t)$ follow from a similar argument as in~\cite[Theorem 24.2]{Rockafellar}.	
If in addition $h$ is continuous, then by applying the fundamental theorem of calculus, one concludes that $\varphi$ is $C^1$ on $(0,\eta)$.\endproof

\begin{proposition}\label{Prop: g-KL modulus is optimal} Let $f:\Rn\rightarrow\overline{\R}$ be proper lsc and let $\bar{x}\in\dom\partial f$. Let $U$ be a nonempty neighborhood of $\bar{x}$ and $\eta\in(0,\infty]$. Let $\varphi\in\Phi_\eta$ be concave and suppose that $f$ has the generalized concave KL property at $\bar{x}$ with respect to $U$, $\eta$ and $\varphi$. Then the exact modulus of the generalized concave KL property of $f$ at $\bar{x}$ with respect to $U$ and $\eta$, denoted by $\tilde{\varphi}$, is well-defined, concave and satisfies \[\tilde{\varphi}(t)\leq\varphi(t),~\forall t\in[0,\eta).\]
Moreover, the function $f$ has the generalized concave KL property at $\bar{x}$ with respect to $U$, $\eta$ and $\tilde{\varphi}$. Consequently, the exact modulus $\tilde{\varphi}$ satisfies
\begin{equation*}
\tilde{\varphi}=\inf\big\{\varphi\in\Phi_\eta:\text{$\varphi$ is a concave desingularizing function of $f$ at $\bar{x}$ with respect to $U$ and $\eta$}\big\}.
\end{equation*}
\end{proposition}
\proof  Let us show first that $\tilde{\varphi}(t)\leq\varphi(t)$ on $[0,\eta)$, which implies immediately that $\tilde{\varphi}$ is well-defined. If $U\cap[0<f-f(\bar{x})<\eta]=\emptyset$, then by our convention $\tilde{\varphi}(t)=0\leq\varphi(t)$ for every $t\in[0,\eta)$. Therefore we proceed with the assumption that $U\cap[0<f-f(\bar{x})<\eta]\neq\emptyset$. By assumption, one has for $x\in U\cap[0<f-f(\bar{x})<\eta]$,
	\[\varphi_-^\prime\big(f(x)-f(\bar{x})\big)\cdot\dist\big(0,\partial f(x)\big)\geq1.\]
which guarantees that $\dist\big(0,\partial f(x)\big)>0$. Fix $s\in(0,\eta)$ and recall from Lemma~\ref{left derivative inequality}(ii) that $\varphi_-^\prime(t)$ is decreasing. Then for $x\in  U\cap[0<f-f(\bar{x})<\eta]$ with $s\leq f(x)-f(\bar{x})$ we have \[\dist^{-1}\big(0,\partial f(x)\big)\leq\varphi_-^\prime\big(f(x)-f(\bar{x})\big)\leq\varphi_-^\prime(s).\]
Taking the supremum over all $x\in U\cap[0<f-f(\bar{x})<\eta]$ satisfying $s\leq f(x)-f(\bar{x})$ yields
	\[h(s)\leq\varphi_-^\prime(s),\]
where $h(s)=\sup\big\{\dist^{-1}\big(0,\partial f(x)\big): x\in U\cap[0<f-f(\bar{x})<\eta],s\leq f(x)-f(\bar{x})\big\}$. If $\lim_{s\rightarrow0^+}h(s)=\infty$, then one needs to treat $\tilde{\varphi}(t)$ as an improper integral. For $t\in(0,\eta)$,
	\[\tilde{\varphi}(t)=\lim_{u\rightarrow0^+}\int_u^th(s)ds\leq\lim_{u\rightarrow0^+}\int_u^t\varphi_-^\prime(s)ds=\varphi(t)<\infty,\]
	where the last equality follows from Lemma~\ref{left derivative inequality}. If $\lim_{s\rightarrow0^+}h(s)<\infty$, then the above argument still applies.
	
 Recall that $\dist\big(0,\partial f(x)\big)>0$ for every $x\in U\cap[0<f-f(\bar{x})<\eta]$. Hence $h(s)$ is positive-valued. Take $s_1,s_2\in(0,\eta)$ with $s_1\leq s_2$. Then for $x\in U\cap[0<f-f(\bar{x})<\eta]$,
\[s_2\leq f(x)-f(\bar{x})\Rightarrow s_1\leq f(x)-f(\bar{x}),\]
implying that $h(s_2)\leq h(s_1)$. Therefore $h(s)$ is decreasing. Invoking Lemma~\ref{Lemma:integration of increasing function is convex}, one concludes that $\tilde{\varphi}$ is a concave function belonging to $\Phi_\eta$, and $\varphi_-^\prime(t)\geq h(t)$ for every $t\in(0,\eta)$.
	
Let $t\in(0,\eta)$. Then for $x\in  U\cap[0<f-f(\bar{x})<\eta]$ with $t=f(x)-f(\bar{x})$,
	\[\tilde{\varphi}_-^\prime\big(f(x)-f(\bar{x})\big)\geq h(t)\geq\dist^{-1}\big(0,\partial f(x)\big),\]
where the last inequality is implied by the definition of $h(s)$, from which the generalized concave KL property readily follows because $t$ is arbitrary.

Recall that $\varphi$ is an arbitrary concave desingularizing function of $f$ at $\bar{x}$ with respect to $U$ and $\eta$, and $\tilde{\varphi}(t)\leq\varphi(t)$ for all $t\in[0,\eta)$. Hence,
\[\tilde{\varphi}\leq\inf\big\{\varphi\in\Phi_\eta:\text{$\varphi$ is a concave desingularizing function of $f$ at $\bar{x}$ with respect to $U$ and $\eta$}\big\}.\]
On the other hand, the converse inequality holds as $\tilde{\varphi}$ is a concave desingularizing function of $f$ at $\bar{x}$ with respect to $U$ and $\eta$.\endproof

Our next example shows that the exact modulus is not necessarily differentiable, which justifies the nonsmooth extension of desingularizing functions in Definition~\ref{Def: g-KL}.

\begin{example}\label{Example: a nonsmooth modulus}Let $\rho>0$. Consider the function given by
		\begin{align*}
			f(x)=
			\begin{cases}
				2\rho|x|-3\rho^2/2,&\text{if }|x|>\rho;\\
				|x|^2/2, &\text{if }|x|\leq\rho.
			\end{cases}
		\end{align*}
		Then the function
		\begin{align*}
			\tilde{\varphi}_1(t)=\begin{cases}\sqrt{2t},&\text{if }0\leq t\leq \rho^2/2;\\ t/(2\rho)+3\rho/4, &\text{if }t>\rho^2/2,\end{cases}
		\end{align*}
		is the exact modulus of the generalized concave KL property of $f$ at $\bar{x}=0$ with respect to $U=\R$ and $\eta=\infty$.
\end{example}
\proof It is easy to see that for $s\in(0,\rho^2/2]$,
\begin{align*}
	h_1(s)&=\sup\big\{\dist^{-1}\big(0,\partial f(x)\big):x\in\R\cap[0<f<\infty],s\leq f(x)\big\}\\
	&=\sup\big\{\dist^{-1}\big(0,\partial f(x)\big):|x|\geq\sqrt{2s}\big\}=1/\sqrt{2s},
\end{align*}
and for $s>\rho^2/2$,
\begin{align*}
	h_1(s)&=\sup\big\{\dist^{-1}\big(0,\partial f(x)\big):x\in\R\cap[0<f<\infty],s\leq f(x)\big\}\\
	&=\sup\big\{\dist^{-1}\big(0,\partial f(x)\big):x\neq0,|x|\geq s/(2\rho)+3\rho/4\big\}=1/(2\rho),
\end{align*}
from which the desired result readily follows.\endproof

It is difficult to compute directly the exact modulus of the generalized concave KL property for multi-variable functions, due to its complicated definition. However, on the real line, we have the following pleasing formula.

\begin{proposition}\label{Prop: modulus of C1 convex} Let $f:\R\rightarrow\overline{\R}$ be proper and lsc. Let $\bar{x}$ be a stationary point. Suppose that there exists an interval $(a,b)\subseteq\inte\dom f$, where $-\infty\leq a<b\leq\infty$, on which $f$ is convex on $(a,b)$ and $C^1$ on $(a,b)\backslash\{\bar{x}\}$. Set $\eta=\min\big\{f(a)-f(\bar{x}),f(b)-f(\bar{x})\big\}$, $f_1(x)=f(x+\bar{x})-f(\bar{x})$ for $x\in(a-\bar{x},0]$ and $f_2(x)=f(x+\bar{x})-f(\bar{x})$ for $[0,b-\bar{x})$. Furthermore, define $\tilde{\varphi}:[0,\eta)\rightarrow\R_+$,
	\begin{equation}\label{Formula:exact modulus for c1convex}
t\mapsto\int_0^t\max\big\{(-f_1^{-1})^\prime(s),(f_2^{-1})^\prime(s)\big\}ds,~\forall t\in(0,\eta)
	\end{equation}	
	and $\tilde{\varphi}(0)=0$. Then $\tilde{\varphi}(t)$ is the exact modulus of the generalized concave KL property at $\bar{x}$ with respect to $U=(a,b)$ and $\eta$. Note that we set $f(x)=\infty$ if $x=\pm\infty$ and $(f_i^{-1})^\prime\equiv0$ if $f_i^{-1}$ does not exist.
\end{proposition}
\proof Replacing $f(x)$ by $g(x)=f(x+\bar{x})-f(\bar{x})$ if necessary, we assume without loss of generality that $\bar{x}=0$ and $f(\bar{x})=0$. Then by the assumption that $\bar{x}=0$ is a stationary point, we have $0\in\partial f(0)=[f_-^\prime(0),f_+^\prime(0)]$, meaning that $f_-^\prime(0)\leq0\leq f_+^\prime(0)$. We learn from
Fact~\ref{Lemma: properties of convex functions on the line} that $f_-^\prime(x)$ and $f_+^\prime(x)$ are increasing functions. Combining the $C^1$ assumption, we have $f^\prime(x)=f_-^\prime(x)\leq f_-^\prime(0)\leq0$ on $(a,0)$ and $f^\prime(x)=f_+^\prime(x)\geq f_+^\prime(0)\geq0$ on $(0,b)$. Hence for $x\in(a,b)\backslash\{0\}$,
\begin{align*}
\dist\big(0,\partial f(x)\big)=|f^\prime(x)|=\begin{cases}-f^\prime(x)=-f_1^\prime(x ),&\text{if }x\in(a,0);\\f^\prime(x)=f_2^\prime(x ),&\text{if }x\in(0,b), \end{cases}
\end{align*}
meaning that the function $x\mapsto\dist\big(0,\partial f(x)\big)$ is decreasing on $(a,0)$ and increasing on $(0,b)$.

Now we work towards showing that $h(s)=\max\big\{(-f_1^{-1})^\prime(s),(f_2^{-1})^\prime(s)\big\}$, where $h(s)$ is the function given in Definition~\ref{Defn: g-KL modulus}. Recall that $f^\prime(x)$ is increasing on $(a,b)\backslash\{0\}$ with $f^\prime(x)\leq0$ on $(a,0)$ and $f^\prime(x)\geq0$ on $(0,b)$. Shrinking the interval $(a,b)$ if necessary, we only need to consider the following four cases.

	\textbf{Case 1:} Suppose that $f_1^\prime(x)<0$ for $x\in(a,0)$ and $f^\prime_2(x)>0$ for $x\in(0,b)$. Then both $f_1$ and $f_2$ are invertible and
	\begin{align*}
	\dist^{-1}\big(0,\partial f(x)\big)=\begin{cases}-1/f_1^\prime(x ),&\text{if }a<x<0;\\ 1/f_2^\prime(x ),&\text{if }0<x<b.\end{cases}
	\end{align*}
	Fix $s\in(0,\eta)$. For $x\in(a,0)$, on which $f_1$ is decreasing,
	\begin{align}\label{ssss}
	s\leq f(x) =f_1(x )\Leftrightarrow f_1^{-1}(s)\geq x .
	\end{align}
	Similarly for $x\in(0,b)$,
	\begin{align}\label{bbbb}
	s\leq f(x) =f_2(x )\Leftrightarrow f_2^{-1}(s)\leq x .
	\end{align}
	Hence one concludes that for $x\in(a,b)$, \[s\leq f(x) \Leftrightarrow x\in(a, f_1^{-1}(s)]\cup[ f_2^{-1}(s),b).\]
On the other hand, we have $0<f(x) <\eta\Leftrightarrow x\in(f_1^{-1}(\eta) ,f_2^{-1}(\eta))\backslash\{0\}$,
where $f_1^{-1}(\eta)>a$ and $f_2^{-1}(\eta)<b$, which means $(a,b)\cap[0<f<\eta]=(f_1^{-1}(\eta),f_2^{-1}(\eta))\backslash\{0\}$.
	
Altogether, we conclude that the function $h:(0,\eta)\rightarrow\R$ given in Definition~\ref{Defn: g-KL modulus} satisfies
	\begin{align*}
	h(s)&=\sup\big\{\dist^{-1}\big(0,\partial f(x)\big):x\in(a,b)\cap[0<f <\eta],s\leq f(x)\big\}\\
	&=\sup\big\{\dist^{-1}\big(0,\partial f(x)\big):x\in( f_1^{-1}(\eta), f_1^{-1}(s)]\cup[ f_2^{-1}(s), f_2^{-1}(\eta))\big\}\\
	&=\max\big\{-1/(f_1^\prime)(f_1^{-1}(s)),1/(f_2^\prime)(f_2^{-1}(s))\big\}\\
	&=\max\big\{(-f_1^{-1})^\prime(s),(f_2^{-1})^\prime(s)\big\},
	\end{align*}
	where the third equality is implied by the fact that
$x\mapsto\dist^{-1}\big(0,\partial f(x)\big)$ is increasing on $(a,0)$ and decreasing on $(0,b)$.
	
	\textbf{Case 2:} If $f^\prime(x)=0$ on $(a,0)$ and $f^\prime(x)>0$ on $(0,b)$, then $f_2$ is invertible and $(f_2^{-1})^\prime(s)=1/f^\prime(f_2^{-1}(s) )>0$ on $(0,\eta)$. Note that by our convention $(f_1^{-1})^\prime(s)$ is set to be zero for all $s$. Hence it suffices to prove $h(s)=(f_2^{-1})^\prime(s)$. For $s\in(0,\eta)$,
	\begin{align*}
	h(s)&=\sup\big\{\dist^{-1}\big(0,\partial f(x)\big):x\in U\cap[0<f <\eta],s\leq f(x)\big\}\\
	&=\sup\big\{1/f_2^\prime(x ): x\in[ f_2^{-1}(s),b)\big\}=1/f_2^\prime(f_2^{-1}(s))=(f_2^{-1})^\prime(s),
	\end{align*}
	where the second equality is implied by (\ref{bbbb}), $U\cap[0<f <\eta]=(0,b)$ and the fact that $1/f_2^\prime(x )$ is decreasing on $(0,b)$.
	
	\textbf{Case 3:} If $f^\prime(x)<0$ on $(a,0)$ and $f^\prime(x)=0$ on $(0,b)$, then $f_1$ is invertible. A similar argument proves that $h(s)=(-f_1^{-1})^\prime(s)$.
	
	\textbf{Case 4:} Now we consider the case where $f^\prime(x)=0$ on $(a,b)$, in which case $U\cap[0<f <\eta]=\emptyset$ and the corresponding exact modulus is $\tilde{\varphi}\equiv0$ by our convention. Moreover, $(-f_1^{-1})^\prime(s)$ and $(f_2^{-1})^\prime(s)$ are set to be constant $0$. Hence we have $\tilde{\varphi}(t)=\int_0^t0~ds=0$, which completes the proof.
\endproof
\begin{remark} In the setting of Proposition~\ref{Prop: modulus of C1 convex}, it is easy to see that the exact modulus satisfies \[\tilde{\varphi}(t)=\int_0^t\frac{ds}{\displaystyle\inf_{(a,b)\cap[f=s]}\dist(0,\partial f(x))},\]
which means that the BDLM desingularizing function given by Fact~\ref{fact: integral condition} coincides with the exact modulus. However, this is not true without the $C^1$ assumption in Proposition~\ref{Prop: modulus of C1 convex}; see Examples~\ref{ex: campare to integral condition with convexity} and~\ref{ex:piewiselinear example}.
\end{remark}
Combining Fact~\ref{Fact: KL at non-stationary point} and Proposition~\ref{Prop: modulus of C1 convex}, we immediately obtain the following corollary.
\begin{corollary} Let $f:(a,b)\rightarrow \R$ be a differentiable convex function. Then $f$ is a concave KL function, i.e., $f$ satisfies
the concave KL-property at every point of $(a,b)$.
\end{corollary}

When proving Proposition~\ref{Prop: modulus of C1 convex}, our initial attempt is to take $\eta>0$ sufficiently small so that $t\mapsto\max\{-f_1^{-1}(t),f_2^{-1}(t)\}$ becomes either $-f_1^{-1}$ or $f_2^{-1}$ on $[0,\eta)$. This attempt leads to a question of independent interest: Let $f$ and $g$ be two smooth strictly increasing convex functions  defined on $[0,\infty)$ with $f(0)=g(0)$.
\[\emph{Is $\inf\{x>0:f(x)=g(x)\}$ always positive?}\]
The answer is negative, as our next example shows.

\begin{example}\label{Prop: inf of intersection}
	There exist strictly increasing convex $C^2$ functions $f,g:\mathbb{R}_+\rightarrow\mathbb{R}_+$ with $f(0)=g(0)=0$ such that $\inf\{x>0:f(x)=g(x)\}=0$. To be specific, let $h:\R\rightarrow\R$ be given by
	 	\begin{align*}
	 	h(x)=
	 	\begin{cases}
	 	\sin\left(\frac{1}{x}\right)e^{-\frac{1}{x^2}},&\text{if }x\neq0;\\
	 	0,&\text{if }x=0.
	 	\end{cases}
	 	\end{align*}
Define $h_+^{\prime\prime}(s)=\max\{h^{\prime\prime}(s),0 \}$ and $h_-^{\prime\prime}(s)=-\min\{h^{\prime\prime}(s),0\}$. Furthermore, set $f_1(x)=\int_0^xh_-^{\prime\prime}(t)dt$ and $g_1(x)=\int_0^xh_+^{\prime\prime}(t)dt$. Then the functions $f,g:\R_+\rightarrow\R_+$ given by
$$g(x)=\int_0^xg_1(t)dt,~f(x)=\int_0^xf_1(t)dt$$
are strictly increasing convex and $C^2$ functions with $f(0)=g(0)=0$, and satisfy $g(x)-f(x)=h(x)$. Hence $\inf\{x>0:f(x)=g(x) \}=\inf\{x>0:h(x)=0\}=0$.
\end{example}

\proof Note that $h(x)\in C^\infty$. We now show that $h$ is a difference of convex functions. Observe from the definition that $h_+^{\prime\prime}(s)$ and $h_-^{\prime\prime}(s)$ are positive-valued and continuous. Then by the fundamental theorem of calculus, $g_1(x)$ and $f_1(x)$ are both increasing $C^1$ functions with $f_1^\prime(x)=h_-^{\prime\prime}(x)$ and $g_1^\prime(x)=h_+^{\prime\prime}(x)$. Since $h^{\prime\prime}(x)=h^{\prime\prime}_+(x)-h^{\prime\prime}_-(x)$, $$g_1(x)-f_1(x)=\int_0^xh^{\prime\prime}(s)ds=h^\prime(x)-h^\prime(0)=h^\prime(x).$$
Suppose that there exists $x_0>0$ such that $g_1(x_0)=0$. Then $g_1(x)=0$ for $x\in(0,x_0)$,
which implies $h^{\prime\prime}(x)\leq0$ on $(0,x_0)$. This is impossible because $h^{\prime\prime}(x)$ oscillates between
positive and negative infinitely many times when $x\rightarrow 0^{+}$.
Hence $g_1$ is strictly positive. A similar argument shows that $f_1$ is also strictly positive. Then applying the fundamental theorem of calculus, one concludes that $f$ and $g$ are strictly increasing $C^2$ functions with $f^\prime(x)=f_1(x)$ and $g^\prime(x)=g_1(x)$. Functions $f$ and $g$ are both convex because $f^{\prime\prime}=f_1^\prime=h^{\prime\prime}_-\geq0$ and $g^{\prime\prime}=g_1^\prime=h^{\prime\prime}_+\geq0$. Furthermore,
	\begin{equation*}
	g(x)-f(x)=\int_0^xh^\prime(s)ds=h(x)-h(0)=h(x).
	\end{equation*}
Hence $\inf\{x>0:f(x)=g(x)\}=\inf\{x>0:h(x)=0\}=\inf\{\frac{1}{n\pi}:n\in\mathbb{N}\}=0$.
\endproof

\subsection{Comparison to the BDLM desingularizing functions}\label{sec: Comparison to integrability condition }
In this subsection, we compare the exact modulus to the BDLM desingularizing functions in Facts~\ref{fact: integral condition} and \ref{fact:integral condition with convexity}. A comparison with the growth condition in~Fact~\ref{Fact:Growth condition gives KL} is also carried out. Examples will be given to show that the exact modulus is the optimal concave desingularizing function, provided that such functions exist.

Below, we use $\varphi$ for the BDLM desingularizing function. By picking a decreasing and continuous majorant\footnote{Note that such majorant may not exist beyond the convex or semialgebraic case.} of $u$ in the integrability condition~(\ref{formula: integral condition}) and integrating, one can get a concave BDLM desingularizing function, which may not be the smallest.

\begin{example} Let $U=\R$ and $\eta=\infty$. Define $f:\R\to\R_+$ by
\begin{align*}
f(x)=
\begin{cases}
\frac{1}{2}x,&\text{if }0\leq x\leq \frac{1}{4};\\
\frac{3}{2}(x-\frac{1}{4})+\frac{1}{8},&\text{if }\frac{1}{4}<x\leq\frac{1}{2};\\
x,&\text{if }x>\frac{1}{2};\\
0&\text{otherwise}.
\end{cases}
\end{align*}
Let $u:(0,\eta)\to\R_+$ be the function given in the integrability condition~(\ref{formula: integral condition}) with $\varepsilon=\bar{r}=\infty$, and let $h$ be given in Definition~\ref{Defn: g-KL modulus}. Then the following statements hold:

(i) Functions $u$ and $h$ are given by
\begin{align*}
	u(s)=
	\begin{cases}
		2,&\text{if }0<s\leq\frac{1}{8};\\
		\frac{2}{3},&\text{if }\frac{1}{8}<s<\frac{1}{2};\\
		1,&\text{if }s\geq\frac{1}{2}.
	\end{cases}\text{ and }h(s)=
\begin{cases}
2,&\text{if }0<s\leq\frac{1}{8};\\
1,&\text{if }s>\frac{1}{8}.
\end{cases}
\end{align*}

(ii) The function $h$ satisfies
$$h=\inf\{\bar{u}:\bar{u}:(0,\eta)\to\R_+\text{ is a continuous and decreasing function with }\bar{u}\geq u\}.$$

(iii) The exact modulus of $f$ at $\bar{x}=0$ with respect to $U$ and $\eta$ is
 \begin{align*}
 	\tilde{\varphi}(t)=
 \begin{cases}
 	2t,&\text{if }0\leq t\leq\frac{1}{8};\\
 	t+\frac{1}{8},&\text{if }t>\frac{1}{8}.
 \end{cases}
 \end{align*}
Let $\bar{u}$ be a continuous and decreasing majorant of $u$ and define $\varphi(t)=\int_0^t\bar{u}(s)ds$. Then $\tilde{\varphi}\leq\varphi$ on $(0,\frac{1}{8}]$ and $\tilde{\varphi}<\varphi$ on $(\frac{1}{8},\infty)$.
\end{example}
\proof (i) The desired results follow from simple calculations. (ii) Define $w=\inf\{\bar{u}:\bar{u}:(0,\eta)\to\R_+\text{ is a continuous and decreasing function with }\bar{u}\geq u\}$. Let $n\in\N$ and define $w_n:(0,\eta)\to\R_+$ by
\begin{align*}
	w_n(s)=
	\begin{cases}
	2,&\text{if }0<s\leq\frac{1}{8};\\
	-n\left(s-\frac{1}{8}\right)+2,&\text{if }\frac{1}{8}<s\leq \frac{1}{8}+\frac{1}{n};\\
	1, &\text{if }s>\frac{1}{8}+\frac{1}{n}.
	\end{cases}
\end{align*}
Then $w_n$ is a decreasing and continuous majorant of $u$ and $w_n\geq w$. Pick $s>\frac{1}{8}$ and note that $\lim_{n\to\infty} w_n(s)=1$. Then $ w(s)\leq\lim_{n\to\infty} w_n(s)=1$, which together with the fact that $w\geq1$ yields $ w(s)=1$. On the other hand, for $s\leq\frac{1}{8}$, we have $2\leq w(s)\leq w_n(s)=2$, which means $w(s)=2$. Therefore, $w=h$ by (i).

(iii) Integrating $h$ yields the desired formula of $\tilde{\varphi}$. Statement(ii) implies that any continuous and decreasing majorant $\bar{u}$ of $u$ satisfies $\bar{u}\geq h$ and there exists some $\varepsilon>0$ such that $\bar{u}(s)>h(s)$ on $(\frac{1}{8},\frac{1}{8}+\varepsilon)$. If there was $t_0>\frac{1}{8}$ such that $\varphi(t_0)=\tilde{\varphi}(t_0)$, then we would have $\bar{u}=h$ almost everywhere on $(0,t_0)$, which is absurd.\endproof

We now compare the exact modulus with Facts~\ref{fact:integral condition with convexity} and \ref{Fact:Growth condition gives KL} by recycling Example~\ref{Example: a nonsmooth modulus}.
On one hand, we shall see that the exact modulus is smaller than any BDLM desingularizing function given by Fact~\ref{fact:integral condition with convexity}. On the other hand, we will show that the smallest desingularizing function obtained from the growth condition in Fact~\ref{Fact:Growth condition gives KL} is still bigger than the exact modulus.

\begin{figure}[H]
	\centering
	\begin{subfigure}{.45\linewidth}
		\includegraphics[width=\textwidth]{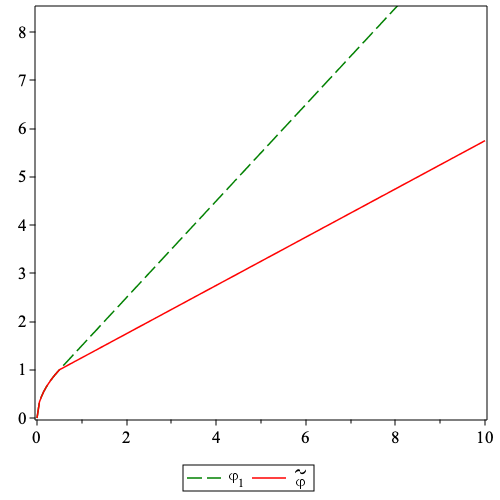}
	\end{subfigure}
	\begin{subfigure}{.45\linewidth}
		\includegraphics[width=\textwidth]{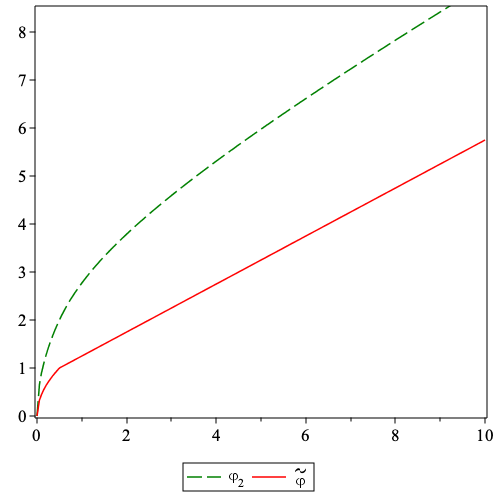}
	\end{subfigure}
	\caption{Plots of Example~\ref{ex: campare to integral condition with convexity}. Left: The exact modulus $\tilde{\varphi}$ and $ \varphi_1$ in the limiting case where $\bar{r}=r_0$. Right: The exact modulus $\tilde{\varphi}$ and $\varphi_2$.}\label{Fig: compare with BDLM}
\end{figure}

\begin{example}\label{ex: campare to integral condition with convexity} Consider the function $f$ given in Example~\ref{Example: a nonsmooth modulus} with $\rho=1$. 	Recall from Example~\ref{Example: a nonsmooth modulus} that the exact modulus of $f$ at $\bar{x}$ with respect to $U=\R$ and $\eta=\infty$ is
	\begin{align*}
		\tilde{\varphi}(t)=\begin{cases}\sqrt{2t},&\text{if }0\leq t\leq 1/2;\\ t/2+3/4, &\text{if }t>1/2.\end{cases}
	\end{align*}
Moreover, the following statements hold:
	
	(i) Applying Fact~\ref{fact:integral condition with convexity} with $r_0=\frac{1}{2}$ and $\bar{r}\in(0,r_0)$ gives that $f$ satisfies the concave KL property at $\bar{x}=0$ with respect to $U=\R$, $\eta=\infty$ and
	\begin{align*}
		 \varphi_1(t)=\begin{cases}\sqrt{2t},&\text{if }t\leq\bar{r};\\\sqrt{2\bar{r}}+\frac{1}{\sqrt{2\bar{r}}}(t-\bar{r}),&\text{if }t>\bar{r}.\end{cases}
	\end{align*}
	Evidently $\tilde{\varphi}\leq \varphi_1$, even in the limiting case where $\bar{r}=r_0$; see the left plot in Figure~\ref{Fig: compare with BDLM}. There are other constructions of $ \varphi_1$, however they are all bigger than the exact modulus $\tilde{\varphi}$, see Remark~\ref{rem:other construction} for a detailed discussion.
	
	(ii) Define $m:\R_+\to\R_+$ by
	\begin{align*}
		m(t)=
		\begin{cases}
			\frac{1}{2}t^2,&\text{if }0\leq t\leq1;\\
			2t-\frac{3}{2},&\text{if }t>1.
		\end{cases}
	\end{align*}
	Then Fact~\ref{Fact:Growth condition gives KL} implies that $f$ has the concave KL property at $0$ with respect to $U=\R$, $\eta=\infty$ and
	\begin{align*}
		\varphi_2(t)=\begin{cases}2\sqrt{2t},&\text{if }t\leq\frac{1}{2};\\\frac{t}{2}+\frac{3}{4}\ln(t)+\frac{7}{4}+\frac{3}{4}\ln(2),&\text{if }t>\frac{1}{2},\end{cases}
	\end{align*}
the smallest one by Fact~\ref{Fact:Growth condition gives KL}. However, $\tilde{\varphi}\leq\varphi_2$; see the right plot in Figure~\ref{Fig: compare with BDLM}.
\end{example}
\proof (i) For $r\in(0,\frac{1}{2}]$, $f(x)=r\Leftrightarrow|x|=\sqrt{2r}$ and for $r\in(\frac{1}{2},\infty)$ we have $f(x)=r\Leftrightarrow|x|=\frac{3}{4}+\frac{1}{2}r$. Then
\begin{align*}
	u(r)=\frac{1}{\displaystyle\inf_{x\in[f=r]}\dist(0,\partial f(x))}=\begin{cases}\frac{1}{\sqrt{2r}},&\text{if }0<r\leq\frac{1}{2};\\\frac{1}{2},&\text{if }r>\frac{1}{2}.\end{cases}
\end{align*}
Noticing that $u$ is continuous on $(0,r_0)$, we set the continuous majorant $\tilde{u}$ in Fact~\ref{fact:integral condition with convexity} to be $u$. The desired $ \varphi_1$ then follows from applying Fact~\ref{fact:integral condition with convexity}.

(ii) Clearly all conditions in Fact~\ref{Fact:Growth condition gives KL} are satisfied. In particular, the equality $f(x)=m(\dist(x,\argmin f))=m(|x|)$ holds for all $x$, which means $m$ is the largest possible modulus of the growth condition. The larger $m$ the smaller its inverse. Hence $\varphi_2(t)=\int_0^t\frac{m^{-1}(s)}{s}ds$ is the smallest possible desingularizing function that one can get from Fact~\ref{Fact:Growth condition gives KL}. The rest of the statement follows from a simple calculation.\endproof
\begin{remark}\label{rem:other construction} The function $\varphi_1$ given in Example~\ref{ex: campare to integral condition with convexity} is indeed $\varphi_1(t)=\int_0^t\bar{u}(s)ds$, where
	\begin{align*}
		\bar{u}(r)=
		\begin{cases}
			u(r),&\text{if }r\leq\bar{r};\\
			u(\bar{r}),&\text{if }r>\bar{r},
		\end{cases}
	\end{align*}
	which is a continuous and decreasing majorant of $u$, where $u$ is given in the proof above. Replacing $\bar{u}$ by other such majorant of $u$ certainly yields a different $\varphi_1$. However, notice that for this example, we have $\tilde{\varphi}(t)=\int_0^t u(s)ds$. Therefore $\tilde{\varphi}(t)\leq\int_0^t\bar{u}(s)ds=\varphi_1(t)$, no matter which majorant $\bar{u}$ we choose.
\end{remark}


Despite the exact modulus $\tilde{\varphi}$ in Example~\ref{ex: campare to integral condition with convexity} is smaller than the desingularizing function obtained from Fact~\ref{fact:integral condition with convexity}, there is still some overlap. In what follows, we construct an example where the exact modulus of a non-differentiable convex function is the strictly smaller one everywhere, except at the origin. Note that we shall compare the exact modulus to Fact~\ref{fact: integral condition} instead of Fact~\ref{fact:integral condition with convexity}, as the former is more general.

\begin{example}\label{ex:piewiselinear example} Let $r_1=\frac{\pi^2}{6}-1$ and $r_{k+1}=r_k-\frac{1}{k^2(k+1)}$ for $k\in\N$. Define for $k\in\N$ and $x>0$
	\[f(x)=\frac{1}{k}\left(x-\frac{1}{k}\right)+r_k,\forall x\in\left(\frac{1}{k+1},\frac{1}{k}\right].\]
Let $f(-x)=f(x)$ for $x<0$ and $f(0)=0$. Then the following statements hold:

(i) The function $f:[-1,1]\rightarrow[0,r_1]$ is continuous and convex with $\argmin f=\{0\}$. 

(ii) The exact modulus of $f$ at $0$ with respect to $U=[-1,1]$ and $\eta=r_1$ is a piecewise linear function $\tilde{\varphi}:[0,r_1]\to\R_+$ satisfying
\begin{align*}	
\tilde{\varphi}(t)=k(t-r_{k+1})+\sum_{i=k+1}^\infty i(r_{i}-r_{i+1}),\forall t\in(r_{k+1},r_k],\forall k\in\N,
\end{align*}
and $\tilde{\varphi}(0)=0$. Furthermore, every desingularizing function $\varphi:[0,r_1]\to\R_+$ obtained from Fact~\ref{fact: integral condition} satisfies $\varphi(t)>\tilde{\varphi}(t)$ on $(0,r_1]$.
\end{example}
\proof(i) Note that we have \[\lim_{k\to\infty}r_{k+1}=r_1-\sum_{i=1}^\infty\frac{1}{i^2(i+1)}=r_1-\sum_{i=1}^\infty\left(\frac{1}{i^2}-\frac{1}{i(i+1)}\right)=r_1-r_1=0,\]
which implies that $f$ is well-defined and continuous at $0$. To see the continuity at $x=\frac{1}{k+1}$ for $k\in\N$, it suffices to observe that
\[\lim_{x\to\frac{1}{k+1}^+}f(x)=\frac{1}{k}\left(\frac{1}{k+1}-\frac{1}{k}\right)+r_k=r_{k+1}=f\left(\frac{1}{1+k}\right),\]
where the second last equality follows from the definition of $r_{k+1}$. Moreover, $f$ is piecewise linear with increasing slope then it is convex.

(ii) Clearly we have $\partial f\left(\frac{1}{k}\right)=\left[\frac{1}{k},\frac{1}{k-1}\right]$ for every $k\in\N$ with $k\geq2$. It follows easily that for $x\in[-1,1]$ with $\frac{1}{k+1}<|x|\leq\frac{1}{k}$, one has
$\dist\left(0,\partial f(x)\right)=\frac{1}{k}$. For $r\in(r_{k+1},r_k]$, $r=f(x)\Leftrightarrow|x|=k(r-r_k)+\frac{1}{k}$. Elementary calculation yields
\begin{align*}
u(r)=\frac{1}{\displaystyle\inf_{x\in U\cap[f=r]}\dist(0,\partial f(x))}=k.
\end{align*}
Note that $u\in L^1(0,r_1)$ hence all conditions in Fact~\ref{fact: integral condition} are satisfied. Indeed,
\[\int_0^{r_1}u(s)ds=\sum_{k=1}^\infty k\cdot(r_k-r_{k+1})=\sum_{k=1}^\infty\frac{1}{k(k+1)} \leq\sum_{k=1}^\infty\frac{1}{k^2}<\infty.\]
Then~\cite[Lemma 44]{bolte2010survey} implies there exists a continuous and decreasing majorant $\bar{u}$ of $u$. The continuity of $\bar{u}$ ensures that for every $k$ there exists $\varepsilon_k>0$ such that $\bar{u}>u$ on $(r_k,r_k+\varepsilon_k)$. Evidently the exact modulus $\tilde{\varphi}$ satisfies $\tilde{\varphi}(t)=\int_0^tu(s)ds$. Altogether, we conclude that the desingularizing function given by Fact~\ref{fact: integral condition} satisfies \[\varphi(t)=\int_0^t\bar{u}(s)ds>\int_0^tu(s)ds=\tilde{\varphi}(t),\forall t\in(0,r_1].\]
Indeed, if there was $t_0\in(0,r_1]$ such that $\varphi(t_0)=\tilde{\varphi}(t_0)$, then we would have $\bar{u}=u$ almost everywhere on $(0,t_0]$, which is absurd.\endproof
\section{The PALM algorithm revisited}\label{Section: PALM}
In this section, we revisit the celebrated \emph{proximal alternating linearized minimization}~(PALM) algorithm. We will show that the exact modulus of the generalized concave KL property leads to the sharpest upper bound on the total length of trajectory of iterates generated by PALM.
\subsection{The PALM algorithm}
Consider the following nonconvex and nonsmooth optimization model:
\begin{equation*}
\min_{(x,y)\in\mathbb{R}^n\times\mathbb{R}^m} \Psi(x,y)=f(x)+g(y)+F(x,y),
\end{equation*}
where $f:\Rn\rightarrow\overline{\R}$ and $g:\R^m\rightarrow\overline{\R}$ are proper and lsc, and $F:\mathbb{R}^n\times\mathbb{R}^m\rightarrow\mathbb{R}$ is $C^1$. This model covers many optimization problems in practice; see~\cite{Bolte2014}. Bolte, Sabach and Teboulle~\cite{Bolte2014} proposed the following algorithm to solve the problem.
\begin{center}
	\fbox{\parbox{\textwidth}
		{\textbf{PALM: Proximal Alternating Linearized Minimization}\\
			1. Initialization: Start with arbitrary $z_0=(x_0,y_0)\in\Rn\times\mathbb{R}^m$.\\
			2. For each $k=0,1,\ldots$, generate a sequence $(z_k)_{k\in\mathbb{N}}=(x_k,y_k)_{k\in\N}$ as follows, where quantities $L_1(y_k)$ and $L_2(x_{k+1})$ will be given in (A2):
			
			$~~$2.1. Take $\gamma_1>1$, set $c_k=\gamma_1L_1(y_k)$ and compute
			\begin{equation}\label{Algorithm: Iteration on x}x_{k+1}\in\prox^f_{c_k}\left(x_k-\frac{1}{c_k}\nabla_xF(x_k,y_k)\right).
			\end{equation}
			
			$~~$2.2. Take $\gamma_2>1$, set $d_k=\gamma_2L_2(x_{k+1})$ and compute
			\begin{equation}\label{Algorithm: Iteration on y}
			y_{k+1}\in\prox^g_{d_k}\left(y_k-\frac{1}{d_k}\nabla_yF(x_{k+1},y_k)\right).
			\end{equation}}
	}
\end{center}

The PALM algorithm is analyzed under the following blanket assumptions in~\cite{Bolte2014}.
\begin{itemize}
	\item[(A1)] $\inf_{\mathbb{R}^n\times\mathbb{R}^m}\Psi>-\infty$, $\inf_{\mathbb{R}^n}f>-\infty$ and $\inf_{\mathbb{R}^m}g>-\infty$.
	
	\item[(A2)]For every fixed $y\in\R^m$, the function $x\mapsto F(x,y)$ is $C^{1,1}_{L_1(y)}$, i.e.,
	\begin{equation*}
	\norm{\nabla_xF(x_1,y)-\nabla_xF(x_2,y)}\leq L_1(y)\norm{x_1-x_2},\forall x_1,x_2\in\mathbb{R}^n.
	\end{equation*}
	Assume similarly that for every $x\in\mathbb{R}^n$, $y\mapsto F(x,y)$ is $C^{1,1}_{L_2(x)}$.
	
	\item[(A3)] For $i=1,2$ there exist $\lambda_i^-,\lambda_i^+>0$ such that
	\begin{align*}
	&\inf\{L_1(y_k):k\in\mathbb{N}\}\geq\lambda_1^-\text{ and }\inf\{L_2(x_k):k\in\mathbb{N}\}\geq\lambda_2^-,\\
	&\sup\{L_1(y_k):k\in\mathbb{N}\}\leq\lambda_1^+\text{ and }\sup\{L_2(x_k):k\in\mathbb{N}\}\leq\lambda_2^+.
	\end{align*}
	\item[(A4)] $\nabla F$ is Lipschitz continuous on bounded subsets of $\mathbb{R}^n\times\mathbb{R}^m$, i.e., on every bounded subset $B_1\times B_2$ of $\mathbb{R}^n\times\mathbb{R}^m$, there exists $M>0$ such that for all $(x_i,y_i)\in B_1\times B_2$, $i=1,2$,
	\begin{equation*}
	\norm{\nabla F(x_1,y_1)-\nabla F(x_2,y_2)} \leq M\norm{(x_1-x_2,y_1-y_2)}.
	\end{equation*}
\end{itemize}
Fact~\ref{fact:prox is well defined} shows that PALM is well defined. Bolte, Sabach and Teboulle showed that the PALM algorithm enjoys the following properties.
\begin{lemma}\label{Lemma: Convergence properties}\emph{\cite[Lemma 3]{Bolte2014}} Suppose that (A1)-(A4) hold. Let $(z_k)_{k\in\N}$ be a sequence generated by PALM. Then the following hold:
	\begin{itemize}
		\item [(i)] The sequence $\big(\Psi(z_k)\big)_{k\in\N}$ is decreasing and in particular
	\begin{align}\label{Formua;function value descent}
	\frac{\rho_1}{2}\norm{z_{k+1}-z_k}^2\leq\Psi(z_k)-\Psi(z_{k+1}),\forall k\geq0,
	\end{align}
where $\rho_1=\min\{(\gamma_1-1)\lambda_1^-,(\gamma_2-1)\lambda_2^-\}$.
		\item[(ii)] $\sum_{k=1}^\infty\norm{z_{k+1}-z_k}^2<\infty$, and hence $\lim_{k\rightarrow\infty}\norm{z_{k+1}-z_k}=0$.
	\end{itemize}
	\end{lemma}
	
	\begin{lemma}\label{Lemma:Approaching critical point}\emph{\cite[Lemma 4]{Bolte2014}} Suppose that (A1)-(A4) hold, and
that $M>0$ is the Lipschitz constant given in~(A4).
Let $(z_k)_{k\in\N}$ be a sequence generated by PALM which is assumed to be bounded. For $k\in\N$, define
		\begin{align*}
		&A_x^k=c_{k-1}(x_{k-1}-x_k)+\nabla_xF(x_k,y_k)-\nabla_xF(x_{k-1},y_{k-1}),\\
		&A_y^k=d_{k-1}(y_{k-1}-y_k)+\nabla_yF(x_k,y_k)-\nabla_yF(x_k,y_{k-1}).
		\end{align*}
	Then $(A_x^k,A_y^k)\in\partial\Psi(x_k,y_k)$, and
	\[\norm{(A_x^k,A_y^k)}\leq\norm{A_x^k}+\norm{A_y^k}\leq(2M+3\rho_2)\norm{z_k-z_{k-1}},\forall k\in\N,\]
	where $\rho_2=\max\{\gamma_1\lambda_1^+,\gamma_2\lambda_2^+\}$.
		\end{lemma}

Denote the set of subsequential limit points of $(z_k)_{k\in\N}$ by $\omega(z_0)=\{z\in\Rn\times\R^m:\exists(z_{k_q})_{q\in\N}\subseteq(z_k)_{k\in\N},z_{k_q}\rightarrow z \text{ as }q\rightarrow\infty \}$. The following lemma summarizes useful properties of $(z_k)_{k\in\mathbb{N}}$ and $\omega(z_0)$, where Lemma~\ref{Lemma: properties of the limit point set}(i) follows from the proof of~\cite[Lemma 5(i)]{Bolte2014}.

		\begin{lemma}\label{Lemma: properties of the limit point set}\emph{\cite[Lemma 5]{Bolte2014}} Suppose that (A1)-(A4) hold. Let $(z_k)_{k\in\N}$ be a sequence generated by PALM which is assumed to be bounded. Then the following assertions hold:
			\begin{itemize}
				\item [(i)] For every $z^*\in\omega(z_0)$ and $(z_{k_q})_{q\in\N}$ converging to $z^*$,  \[\lim_{q\rightarrow\infty}\Psi(z_{k_q})=\Psi(z^*).\]
				Moreover, $\omega(z_0)\subseteq\text{stat }\Psi$, where $\text{stat }\Psi$ denotes the set of stationary points of $\Psi$.
				\item[(ii)] $\lim_{k\rightarrow\infty}\dist(z_k,\omega(z_0))=0$.
				\item[(iii)] The set $\omega(z_0)$ is nonempty, compact and connected.
				\item[(iv)] The objective function is constant on $\omega(z_0)$.
			\end{itemize}
			\end{lemma}
\subsection{The sharpest upper bound for the total length of trajectory of iterates}

In this subsection, we improve a result by Bolte, Sabach and Teboulle~\cite[Theorem 1]{Bolte2014}. We begin with a technical lemma, which is a sharper version of \cite[Lemma 6]{Bolte2014}.
\begin{lemma}\label{Lemma: Uniformized g-KL and modulus} Let $f:\Rn\rightarrow\overline{\R}$ be proper lsc and let $\mu\in\R$. Let $\Omega\subseteq\dom\partial f$ be a nonempty compact set on which $f(x)=\mu$ for all $x\in\Omega$. Suppose that $f$ has the pointwise generalized concave KL property at each $x\in\Omega$. Let $\varepsilon$, $\eta>0$ and concave $\varphi\in\Phi_\eta$ be those given in Proposition~\ref{Prop: Uniformize the g-KL}. Set $U=\Omega_\varepsilon$ and define $h:(0,\eta)\rightarrow\R_+$ by \[h(s)=\sup\big\{\dist^{-1}\big(0,\partial f(x)\big):x\in U\cap[0<f-\mu<\eta],s\leq f(x)-\mu\big\}.\]
Then the function $\tilde{\varphi}:[0,\eta)\rightarrow\R_+ :t\mapsto\int_0^th(s)ds,~\forall t\in(0,\eta)$,
with $\tilde{\varphi}(0)=0$, is well-defined, concave and belongs to $\Phi_\eta$. The function $f$ has the setwise generalized concave KL property on $\Omega$ with respect to $U$, $\eta$ and $\tilde{\varphi}$. Consequently,
\[\tilde{\varphi}=\inf\big\{\varphi\in\Phi_\eta:\text{$\varphi$ is a concave desingularizing function of $f$ on $\Omega$ with respect to $U$ and $\eta$}\big\}.\]
We say $\tilde{\varphi}$ is the exact modulus of the setwise generalized concave KL property of $f$ on $\Omega$ with respect to $U$ and $\eta$.
\end{lemma}
\proof Apply a similar argument as in Proposition~\ref{Prop: g-KL modulus is optimal}. \endproof

The following theorem provides the ``sharpest" upper bound for the total length of the trajectory of iterates generated by PALM, which improves Bolte, Sabach and Teboulle~\cite[Theorem 1]{Bolte2014}. The notion of ``sharpest" will be specified later in Remark~\ref{rem: sharpest}. Our proof follows a similar approach as in~\cite[Theorem 1]{Bolte2014}, but makes use of the exact modulus of the setwise generalized concave KL property.

\begin{theorem}\label{Theorem: finite length property} Suppose that the objective function $\Psi$ is a generalized concave KL function such that (A1)-(A4) hold. Let $(z_k)_{k\in\N}$ be a sequence generated by PALM which is assumed to be bounded. Then the following assertions hold:
	\begin{itemize}
\item[(i)] The sequence $(z_k)_{k\in\N}$ converges to a stationary point $z^*$ of objective function $\Psi$.
		\item[(ii)] The sequence $(z_k)_{k\in\N}$ has finite length. More precisely, there exist
$l\in\N$, $\eta\in(0,\infty]$ and $\tilde{\varphi}\in\Phi_\eta$ such that for $p\geq l+1$ and every $q\in\N$
		\begin{equation}\label{Formula:upper bound of partial sum}
		\sum_{k=p}^{p+q}\norm{z_{k+1}-z_k}\leq C\cdot\tilde{\varphi}\big(\Psi(z_p)-\Psi(z^*)\big)+\norm{z_p-z_{p-1}}.
		\end{equation}
		Therefore
		\begin{equation}\label{Formula:Sharpest upper bound}\sum_{k=1}^\infty\norm{z_{k+1}-z_k}\leq A+C\cdot\tilde{\varphi}\big(\Psi(z_{l+1})-\Psi(z^*)\big)<\infty,
		\end{equation}
		where $A=\norm{z_{l+1}-z_{l}}+\sum_{k=1}^{l}\norm{z_{k+1}-z_k}<\infty$ and $C=2(2M+3\rho_2)/\rho_1$.
\end{itemize}
\end{theorem}
\proof Because $(z_{k})_{k\in \N}$ is bounded, there exists a convergent subsequence,
say $z_{k_q}\rightarrow z^*\in\omega(z_0)$. Then Lemma~\ref{Lemma: properties of the limit point set}(i) implies $\lim_{q\rightarrow\infty}\Psi(z_{k_q})=\Psi(z^*)$ and $z^*\in \text{crit }\Psi$.
Since $(\Psi(z_k))_{k\in\N}$ is a decreasing sequence
by Lemma~\ref{Lemma: Convergence properties},
we have $\lim_{k\rightarrow\infty}\Psi(z_k)=\lim_{q\rightarrow\infty}
\Psi(z_{k_q})=\Psi(z^*)$.

We will show that $(z_{k})_{k\in \N}$ converges to $z^*$, and along the way we also establish
\eqref{Formula:upper bound of partial sum} and \eqref{Formula:Sharpest upper bound}.
We proceed by considering two cases.

\textbf{Case 1:} If there exits $l$ such that $\Psi(z_{l})=\Psi(z^*)$, then by the decreasing property of $(\Psi(z_k))_{k\in\N}$, one has $\Psi(z_{l+1})=\Psi(z_{l})$ and therefore $z_{l}=z_{l+1}$ by~(\ref{Formua;function value descent}). Hence by induction, we conclude that $\lim_{k\rightarrow\infty}z_k=z^*$. The desired assertion follows immediately.

\textbf{Case 2:} Now we consider the case where $\Psi(z^*)<\Psi(z_k)$ for all $k\in\N$. By Lemma~\ref{Lemma: properties of the limit point set} and assumption, $\Psi$ is a generalized concave KL function that is constant on compact set $\omega(z_0)$. Invoking Lemma~\ref{Lemma: Uniformized g-KL and modulus} shows that there exist $\varepsilon>0$ and $\eta>0$ such that the exact modulus of the setwise generalized concave KL property on $\Omega=\omega(z_0)$ with respect to $U=\Omega_\varepsilon$ and $\eta$ exists, which is denoted by $\tilde{\varphi}$. Hence for every $z\in \Omega_\varepsilon\cap[0<\Psi-\Psi(z^*)<\eta]$,
\begin{align}\label{Formula: use g-KL in finite length-0}
\tilde{\varphi}_-^\prime\big(\Psi(z)-\Psi(z^*)\big)\cdot\dist\big(0,\partial\Psi(z)\big)\geq1.
\end{align}
Since
 $\lim_{k\rightarrow\infty}\Psi(z_k)=\Psi(z^*)$, there exists some $l_1>0$ such that $0<\Psi(z_k)-\Psi(z^*)<\eta$ for $k>l_1$. On the other hand, Lemma~\ref{Lemma: properties of the limit point set}(ii) shows that there exists $l_2>0$ such that $\dist(z_k,w(z_0))<\varepsilon$ for $k>l_2$. Altogether, we conclude that for
 $k>l=\max\{l_1,l_2\}$, $z_k\in\Omega_\varepsilon\cap[0<\Psi-\Psi(z^*)<\eta]$ and
\begin{equation}\label{Formula: use g-KL in finite length}
\tilde{\varphi}_-^\prime\big(\Psi(z_k)-\Psi(z^*)\big)\cdot\dist\big(0,\partial\Psi(z_k)\big)\geq1.
\end{equation}
It follows from Lemma~\ref{Lemma:Approaching critical point} that $\dist\big(0,\partial\Psi(z_k)\big)\leq\norm{\left(A_x^k,A_y^k\right)}\leq(2M+3\rho_2)\norm{z_k-z_{k-1}}$. Hence one has from (\ref{Formula: use g-KL in finite length}) that for $k>l$,
\begin{equation}\label{Formula:left derivative lower bound}
\tilde{\varphi}^\prime_-\big(\Psi(z_k)-\Psi(z^*)\big)\geq\dist^{-1}\big(0,\partial\Psi(z_k)\big)\geq\frac{1}{2M+3\rho_2}\norm{z_k-z_{k-1}}^{-1}.
\end{equation}
Note that $\norm{z_k-z_{k-1}}\neq0$. Otherwise Lemma~\ref{Lemma:Approaching critical point} would imply that $$\dist\big(0,\partial\Psi(z_k)\big)\leq(2M+3\rho_2)\norm{z_k-z_{k-1}}=0,$$ which contradicts to (\ref{Formula: use g-KL in finite length}). Applying Lemma~\ref{left derivative inequality}(ii) to $\tilde{\varphi}$ with $s=\Psi(z_{k+1})-\Psi(z^*)$ and $t=\Psi(z_k)-\Psi(z^*)$, one obtains for $k>l$
\begin{align}\label{Formula:left derivative upper bound}
\nonumber\frac{\tilde{\varphi}\big(\Psi(z_k)-\Psi(z^*)\big)-\tilde{\varphi}\big(\Psi(z_{k+1})-\Psi(z^*)\big)}{\Psi(z_k)-\Psi(z_{k+1})}&\geq\tilde{\varphi}^\prime_-\big(\Psi(z^{k})-\Psi(z^*)\big)\\
&\geq\frac{1}{2M+3\rho_2}\norm{z_k-z_{k-1}}^{-1}.
\end{align}
For the sake of simplicity, we set
\begin{equation*}
\Delta_{p,q}=\tilde{\varphi}\big(\Psi(z_p)-\Psi(z^*)\big)-\tilde{\varphi}\big(\Psi(z_q)-\Psi(z^*)\big).
\end{equation*}
Then (\ref{Formula:left derivative upper bound}) can be rewritten as
\begin{equation}
\Psi(z_k)-\Psi(z_{k+1})\leq\norm{z_k-z_{k-1}}\cdot\Delta_{k,k+1}\cdot(2M+3\rho_2).
\end{equation}
Furthermore, Lemma~\ref{Lemma: Convergence properties}(i) gives
\begin{align}\label{Formula:upper bound for norm{z_k+1-z_k}^2}
\norm{z_{k+1}-z_k}^2\leq\frac{2}{\rho_1}\left[\Psi(z_k)-\Psi(z_{k+1})\right]\leq C\Delta_{k,k+1}\norm{z_k-z_{k-1}},
\end{align}
where $C=\frac{2(2M+3\rho_2)}{\rho_1}\in(0,\infty)$. By the geometric mean inequality $2\sqrt{\alpha\beta}\leq\alpha+\beta$ for $\alpha,\beta\geq0$, one gets for $k>l$
\begin{align}
2\norm{z_{k+1}-z_k}\leq C\Delta_{k,k+1}+\norm{z_k-z_{k-1}}.
\end{align}

Let $p\geq l+1$. For every $q\in\mathbb{N}$, summing up the above inequality from $p$ to $p+q$ yields
\begin{align*}
2\sum_{k=p}^{p+q}\norm{z_{k+1}-z_k}&\leq C\sum_{k=p}^{p+q}\Delta_{k,k+1}+\sum_{k=p}^{p+q}\norm{z_k-z_{k-1}}+\norm{z_{{p+q}+1}-z_{p+q}}\\
&=C\Delta_{p,p+q+1}+\sum_{k=p}^{p+q}\norm{z_{k+1}-z_k}+\norm{z_p-z_{p-1}}\\
&\leq C\tilde{\varphi}\big(\Psi(z_p)-\Psi(z^*)\big)+\sum_{k=p}^{p+q}\norm{z_{k+1}-z_k}+\norm{z_p-z_{p-1}},
\end{align*}
where the last inequality holds because $\tilde{\varphi}\geq0$. Hence for $q\in\N$
\begin{equation*}
\sum_{k=p}^{p+q}\norm{z_{k+1}-z_k}\leq C\cdot\tilde{\varphi}\big(\Psi(z_p)-\Psi(z^*)\big)+\norm{z_p-z_{p-1}},
\end{equation*}
which proves (\ref{Formula:upper bound of partial sum}). By taking $q\rightarrow\infty$,
\begin{equation*}\sum_{k=1}^\infty\norm{z_{k+1}-z_k}\leq\sum_{k=1}^{p-1}\norm{z_{k+1}-z_k}+C\cdot\tilde{\varphi}\big(\Psi(z_p)-\Psi(z^*)\big)+\norm{z_p-z_{p-1}},
\end{equation*}
from which (\ref{Formula:Sharpest upper bound}) readily follows by setting $p=l+1$.

Now let $p\geq l+1$, where $l$ is the index given in assertion (i), and let $q\in\N$. Then
\[\norm{z_{p+q}-z_{p}}\leq\sum_{k=p}^{p+q-1}\norm{z_{k+1}-z_k}\leq\sum_{k=p}^{p+q}\norm{z_{k+1}-z_k}.\]
Recall that $\tilde{\varphi}(t)\rightarrow0$ as $t\rightarrow0^+$, $\Psi(z_k)-\Psi(z^*)\rightarrow0$ and $\norm{z_{k+1}-z_k}\rightarrow0$ as $k\rightarrow\infty$. Invoking (\ref{Formula:upper bound of partial sum}), one obtains that
\[\norm{z_{p+q}-z_{p}}\leq C\cdot\tilde{\varphi}\big(\Psi(z_p)-\Psi(z^*)\big)+\norm{z_p-z_{p-1}}\rightarrow0,p\rightarrow\infty,\]
meaning that $(z_k)_{k\in\N}$ is Cauchy and hence convergent.
Because $z_{k_{q}}\rightarrow z^*$, we conclude that $z_{k}\rightarrow z^*$.
\endproof

\begin{remark}\label{rem: sharpest} The bound for the total length of iterates (\ref{Formula:Sharpest upper bound}) is the ``sharpest", in the sense that it is the smallest one can get by using the usual KL convergence analysis. Assuming that the objective function $\Psi$ has the concave KL property, Bolte, Sabach and Teboulle~\cite[Theorem 1]{Bolte2014} showed that \begin{equation}\label{upper bound by bolte}\sum_{k=1}^\infty\norm{z_{k+1}-z_k}\leq A+C\cdot\varphi\big(\Psi(z_{l+1})-\Psi(z^*)\big),\end{equation}
where $\varphi(t)$ is in our terminology a concave desingularizing function for the setwise concave KL property of $\Psi$ on $\Omega=w(z_0)$ with respect to $\Omega_\varepsilon$ and $\eta>0$; see~\cite[Lemma 6]{Bolte2014}. Note that $A$ and $C$ are fixed. Then we learn from~(\ref{upper bound by bolte}) that the smaller $\varphi(t)$ is, the sharper the upper bound becomes. According to Lemma~\ref{Lemma: Uniformized g-KL and modulus}, $\tilde{\varphi}$ is the smallest among all possible $\varphi$. Hence the upper bound given by~(\ref{Formula:Sharpest upper bound}) is the sharpest.
\end{remark}

\section{Conclusion}\label{Sec:conclusion}

In this work, we introduced the generalized concave KL property and its exact modulus, which answers the open question~(\ref{open question}). Our results open the door for obtaining sharp results of algorithms that adopt the concave KL assumption. We conclude this paper with some future directions:
\begin{itemize}	

\item Compute or at least estimate the exact modulus of the generalized concave KL property for concrete optimization models.
\item One way to estimate the exact modulus is applying calculus rules of the generalized concave KL property. Li and Pong~\cite{li2018calculus} and Yu et~al.~\cite{yu2019deducing} developed several calculus rules of the concave KL property, in the case where desingularizing functions take the specific form $\varphi(t)=c\cdot t^{1-\theta}$, where $c>0$ and $\theta\in[0,1)$. However, the exact modulus has various forms depending on the given function, which requires us to obtain general calculus rules without assuming desingularizing functions of any specific form.
\end{itemize}

\section*{Acknowledgments}
XW and ZW were partially supported by NSERC Discovery Grants. The authors thank Dr.\ Heinz H. Bauschke for many useful
discussions and suggesting this research problem. They are also very grateful to editor Michael P.\ Friedlander and
anonymous referees for their valuable comments and suggestions on terminologies that improved this manuscript significantly.

\bibliographystyle{siam}
\bibliography{KL_modulus_reference}

\end{document}